\documentclass[preprint,11pt]{elsarticle}

\title{A Hybrid High-Order Finite Element Method for a Nonlocal Nonlinear Problem of Kirchhoff Type}

 \usepackage{amsmath,amsthm,amssymb,enumerate}
\usepackage{gensymb}

\usepackage[sc]{mathpazo}
\usepackage{multirow} 

\usepackage{multicol}

\usepackage{newtxtext,newtxmath}
\usepackage{subfig}
\usepackage{graphicx}
\usepackage{epstopdf}
\usepackage{hyperref}
\usepackage{cancel}
\usepackage[margin=3cm]{geometry}
\usepackage{tikz}
\usepackage{caption}
\usetikzlibrary{shapes,calc}
\usepackage{verbatim}
\usepackage{mathrsfs}
\usepackage{algorithm}
\usepackage{accents}

\newtheorem{theorem}{Theorem}[section]
\newtheorem{lemma}[theorem]{Lemma}
\newtheorem{proof of lemma}[theorem]{Proof of Lemma}
\newtheorem{proposition}[theorem]{Proposition}

\theoremstyle{definition}

\newtheorem{example}[theorem]{Example}

\newtheorem{remark}[theorem]{Remark}
\numberwithin{equation}{section}

\newcommand{\dx}{{\rm\,dx}}

\newcommand{\bfnTF}{\boldsymbol{n}_{TF}}
\newcommand{\bfn}{\boldsymbol{n}}

\newcommand{\mR}{\mathbb R}

\newcommand{\cF}{\mathcal F}

\newcommand{\cT}{\mathcal T}

\newcommand{\cM}{\mathcal M}

\newcommand{\fl}{\quad \forall\:}

\begin{document}
\begin{frontmatter}
\author{Gouranga Mallik\footnote{Department of Mathematics, School of Advanced Sciences, Vellore Institute of Technology, Vellore - 632014, India. Email. gouranga.mallik@vit.ac.in.}}
\begin{abstract}
In this article, we design and analyze a hybrid high-order (HHO) finite element approximation for the solution of a nonlocal nonlinear problem of Kirchhoff type. The HHO method involves arbitrary-order polynomial approximations on structured and unstructured polytopal meshes. We establish the existence of a unique discrete solution to the nonlocal nonlinear discrete problem. We derive an optimal-order error estimate in the discrete energy norm. The discrete system is solved using Newton's iterations on the sparse matrix system. We perform numerical tests to substantiate the theoretical results.
\end{abstract}
\end{frontmatter}

\section{Introduction}\label{intro_Goal_bih}
\subsection{Model problem}
We consider a nonlocal nonlinear elliptic boundary value problem of Kirchhoff type: 
\begin{align}\label{cts_kirchhoff}
    -M\left(\int_{\Omega}|\nabla u|^2\dx\right)\Delta u = f(x,u)\;\;\text{ in }\; \Omega,\;\; u=0\;\;\text{ in }\;\partial\Omega,
\end{align}
where $\Omega\subset \mR^d,\, d=2,3$ is a bounded polytopal domain with Lipschitz boundary $\partial\Omega$, the Euclidean norm of the gradient, $|\nabla u|^2=(\frac{\partial u}{\partial x})^2+(\frac{\partial u}{\partial y})^2$, the Laplacian operator $\Delta u=\frac{\partial^2u}{\partial x^2}+\frac{\partial^2u}{\partial y^2}$ and the given source term $f:\Omega\times\mR\to\mR$ and the nonlocal coefficient term $M:\mR_{+}\to\mR$. Equation \eqref{cts_kirchhoff} involves an integral over $\Omega$, implying that the problem is a nonlocal and does not hold pointwise. The main aim of the article is to design and analyze hybrid high-order polytopal finite element approximation for the Kirchhoff type problem~\eqref{cts_kirchhoff}. We establish the existence and uniqueness of the discrete problem and derive an error estimate.

\subsection{Literature survey}
The above problem~\eqref{cts_kirchhoff} corresponds to the stationary form of the generalized Kirchhoff equation:
\begin{equation}\label{cts_wave_Kirch}
    \frac{\partial^2u}{\partial t^2}-M\left(\int_{\Omega}|\nabla u|^2\dx\right)\Delta u = f(x,u),
\end{equation}
which is an extension of the classical wave equation that accounts for volume changes in a vibrating body. In \cite{Chipot92_Nonlocal,Gudi12_Kirchhof_apriori}, problems structurally related to \eqref{cts_kirchhoff}, where the nonlocal term is of the form $M(\int_\Omega u \dx)$, model certain biological phenomena of the density of bacteria in the domain $\Omega$. The well-posedness and qualitative properties of the problems \eqref{cts_kirchhoff} and \eqref{cts_wave_Kirch} are studied in \cite{Dancona92_Solvability_Kirchhoff,Arosio_96_well-posedness_Kirchhoff}; see also the survey article \cite{MA05_Kirchhoff}.
During the past two decades, the numerical approximation of Kirchhoff-type problems of the forms~\eqref{cts_kirchhoff} and~\eqref{cts_wave_Kirch} has become a growing area of research. In \cite{Peradze_05_Spectral_Kirchhoff}, Peradze analyzed the one-dimensional Kirchhoff string problem using a spectral method. A seminal contribution by Gudi~\cite{Gudi12_Kirchhof_apriori} to the finite element approximation of the Kirchhoff type problem~\eqref{cts_kirchhoff} involves a reformulation into an equivalent discrete system, in which the Jacobian arising from the linearization preserves sparsity. The author established the existence and uniqueness of the discrete solution, together with optimal-order error estimates in the energy norm. In \cite{DondPani17_Kirchhoff}, the authors discuss a priori and a posteriori error estimates of conforming and mixed FEMs for the nonlocal Kirchhoff problem of type \eqref{particular_kirchhoff}. In \cite{Adak_Natarajan20_VEM_Kirchhoff}, a virtual element method for \eqref{cts_kirchhoff} is considered in a global $H^1$-conforming virtual element space. Recently, in \cite{Gudi_Pal25_AFEM_Kirchhof}, the authors consider an a posteriori error estimate for \eqref{cts_kirchhoff} and prove the convergence of the adaptive FEM.

Over the past decades, significant progress has been made in developing high-order polynomial approximation techniques for partial differential equations on general polytopal meshes. Among these, several well-established approaches have been developed and refined. The Hybridizable Discontinuous Galerkin (HDG) method was proposed in \cite{Cock_Gopal_Laz_09_Unif}. The Virtual Element Method (VEM) has been extensively investigated in \cite{Beira_Brezz_Cang_Man_Mar_13_VEM,Beira_Brezz_Mari_13_VEM_Elast,Brezz_Falk_Mari_14_VEM_Mixed}.
Other notable frameworks include the Weak Galerkin method \cite{Wang_Ye_13_Weak_Galerkin,Wang_Ye_14_Mixed_Weak_Galerkin,Mu_Wang_Ye_15_WeakGaler_Nonlin} and the Gradient Discretization method \cite{Dron_Eym_Herb_GDM_16,Dron_Eym_GDM_18_book,Piet_Dron_Man_18_GradDis_Polytope}, and the Multiscale Hybrid-Mixed method \cite{Arya_Hard_Pare_Vale_13_Multiscale}. The present work focuses on the Hybrid High-Order (HHO) method, first introduced in \cite{Piet_Ern_Lem_14_arb_local,Di_Pietr_Ern_HHO_el_15} and reviewed systemically in \cite{Piet_Jero_HHO_Book_20}. The HHO framework offers a dimension-independent design on polytopal meshes based on local polynomial reconstruction and supports element-wise static condensation, enabling significant reductions in global system size and computational cost. Because of these properties, HHO method has emerged as a versatile and efficient technique for the numerical approximation of PDEs on complex meshes.

The Hybrid High-Order (HHO) method shares several conceptual similarities with the Hybridizable Discontinuous Galerkin (HDG) method, yet the two differ primarily in their choice of stabilization; see \cite{Cock_Piet_Ern_16_HHO_HDG} for further details. In the nonconforming Virtual Element Method (ncVEM), stabilization is based on the projection of the virtual function, whereas in the HHO framework, stabilization relies on a reconstructed polynomial representation of the function. The relationship between the HHO and VEM approaches has been analyzed in \cite{Lem_21_Bridging_HHO_VEM}. In the lowest-order case ($k=0$), the HHO method is closely aligned with the Hybrid Mixed Mimetic family, including the Mimetic Finite Difference, Hybrid Finite Volume, and Mixed Finite Volume methods; see \cite{Dron_Eym_Gall_10_unified, Brez_Lipn_Simo_05_Mim_Dif,Eym_Gall_Herb_10_NonconfMesh,Dron_14_FV_Diff, Dron_Eym_06_MixedFV}. Significant contributions to the development of HHO methods include their applications to a variety of problems: pure diffusion \cite{Piet_Ern_Lem_14_arb_local}, linear elasticity \cite{Di_Pietr_Ern_HHO_el_15}, interface problems \cite{Burman_Ern_18_HHO_Interface}, elliptic obstacle problem \cite{Cicu_Ern_Gudi_20_HHO_Obstacle} and the Stokes system \cite{Piet_Ern_Link_16_HHO_Stokes}.

The Hybrid High-Order (HHO) method has also been successfully applied to various classes of nonlinear problems, including the steady Leray–Lions equations~\cite{Piet_Dron_16_Leray}, the Oseen problem~\cite{Aghili18_Oseen_HHO,GM_RB_TG_24_Oseen_LPS}, and the steady incompressible Navier–Stokes equations~\cite{Piet_Krell_18_NSE}. Further developments have extended the method to quasilinear elliptic problems of nonmonotone type~\cite{TG_GM_TP_22_HHO_Quasi} and to strongly nonlinear elliptic models~\cite{GM_TG_24_HHO_Strong_Nonlin}. In this paper, we apply the Hybrid High-Order framework to nonlocal nonlinear problems of Kirchhoff type.

\subsection{Contributions}
The contributions of this work are summarized as follows:
\begin{itemize}
    \item \textbf{Design of a polytopal HHO method for Kirchhoff type problems:} We propose a hybrid high-order polytopal finite element method for the nonlocal nonlinear problem of Kirchhoff type. The method offers flexibility in choosing arbitrary order polynomial approximations on polytopal meshes. The HHO method is based on the local polynomial spaces in cells and on faces to achieve high-order accuracy while preserving computational efficiency through static condensation.
    \item \textbf{Existence, uniqueness, and error estimate for the discrete HHO approximation:} The existence and uniqueness of the discrete HHO solution are established without additional assumptions on the continuous problem. Furthermore, for the particular nonlocal Kirchhoff problem~\eqref{particular_kirchhoff}, we prove the existence and uniqueness of the discrete solution for a modified HHO formulation without requiring any assumptions on the continuous problem. We derive an a priori error estimate in the discrete energy norm. The analysis shows optimal convergence rates with respect to the polynomial degree under sufficient regularity of the solution, thus extending the theoretical framework of HHO methods to a class of nonlocal nonlinear problems. The analysis highlights that for the case of (nonconforming) HHO schemes, a sufficiently small mesh parameter is required to control the nonlocal and nonlinear features of the Kirchhoff Type problem and to handle the inconsistency of the discrete formulation.
    \item \textbf{Numerical scheme and tests:} We describe an equivalent numerical scheme in the context of the HHO method that ensures the sparse matrix structure of the underlying Jacobian matrix to solve the discrete problem. Numerical tests are performed on polygonal meshes (triangular, Cartesian, hexagonal, and Kershaw meshes) with polynomial degrees $k=0,1,2$ to illustrate the convergence rate and robustness of the present method.
\end{itemize}

\subsection{Structure of the paper}
The article is organized as follows. Section~\ref{sect:cts} discusses the existence and uniqueness of the solution to the continuous problem~\eqref{cts_kirchhoff} under standard assumptions. Preliminaries on the HHO method are provided in Section~\ref{sect:HHO_Dis}. In Section~\ref{sect:dis_HHO_form}, we design the HHO method for the Kirchhoff problem~\eqref{cts_kirchhoff} and prove the existence of a discrete unique solution. The optimal order error estimate is derived in Section~\ref{sect:err_est}. Finally, in Section~\ref{sect:numer_exp}, numerical experiments are performed to validate the theoretical results.

Throughout the paper, we employ standard notation for Lebesgue and Sobolev spaces and their associated norms. For a subdomain $K \subset \Omega$, the $L^2$-inner product on $L^2(K)$ is denoted by $(\cdot, \cdot)_K$, and the corresponding $L^2$-norm by $\|\cdot\|_K$. When $K = \Omega$, the subscript is omitted. The usual seminorm and norm on $H^{s}(\Omega)$ for $s > 0$ are denoted by $|\cdot|_{s}$ and $\|\cdot\|_{s}$. 

%\begin{itemize}
%\item This paper designs a hybrid high-order polytopal finite element method for the nonlocal nonlinear problem of Kirchhoff type.
% \item The existence and uniqueness of the discrete HHO solution are established without additional assumptions about the continuous problem. Sufficiently small mesh parameter is imposed to handle inconsistency of the HHO method like the case of conforming finite element methods.
%\item A prior error estimate in the discrete energy norm is derived.
% \item Numerical experiments are performed to validate the theoretical results.
% \end{itemize}

\section{Continuous problem and week formulation}\label{sect:cts}
We make the following assumptions \cite{Gudi12_Kirchhof_apriori,Gudi_Pal25_AFEM_Kirchhof} on the data which are essential for the existence and uniqueness of the solution of \eqref{cts_kirchhoff}:

\noindent A1: The source function $f:\Omega\times\mR\to\mR$ is a Carath\'{e}odory function and Lipschitz continuous for some constant $L_f>0$
\begin{equation*}
    |f(x,s)-f(x,t)|\leq L_f |s-t|\quad\forall x\in\Omega, \quad\forall s,t\in\mR.
\end{equation*}

\noindent A2: The nonlocal function $M:\mR_{+}\to\mR$ is continuous and there exists a constant $m_0>0$ such that
\begin{equation*}
    M(s)\geq m_0\quad\forall s\geq 0,
\end{equation*}
and $m_0$ is such that $m_0\lambda_1> L_f$, where $\lambda_1$ is the principal eigenvalue of the Laplace operator $-\Delta$ over $H^1_0(\Omega)$.

\noindent A3: Let 
\begin{equation}\label{def_R1}
   R_1=\frac{\|f(\cdot,0)\|}{\lambda_1^{1/2}(m_0-L_f\lambda_1^{-1})}.
\end{equation}
Assume $M$ is a Lipschitz continuous function in $[0,R_1^2]$ with the Lipschitz constant $L_M$, and $L_M$ is such that $(m_0-L_f\lambda_1^{-1}-2L_MR_1^2)>0$.

A weak formulation of \eqref{cts_kirchhoff} is given by: Find $u\in H^1_0(\Omega)$ such that 
\begin{equation}\label{weak_kirchhoff}
  -M\left(\int_{\Omega}|\nabla u|^2\dx\right)(\nabla u,\nabla v)=(f(u),v)\quad\forall v\in H^1_0(\Omega).
\end{equation}
Here we interpret $(f(u),v)$ as $\int_{\Omega}f(x,u)v(x)\dx$. The following theorem on the existence and uniqueness of the solution to \eqref{cts_kirchhoff} follows from \cite[Theorem 1 \& 2]{MA05_Kirchhoff}.

\begin{theorem}
\cite{Gudi12_Kirchhof_apriori} Assume that (A1)--(A2) hold. Then problem \eqref{weak_kirchhoff} has at least one solution in $H^1_0(\Omega)$. Moreover, any solution of \eqref{weak_kirchhoff} satisfies
\begin{equation*}
    \|\nabla u\|\leq R_1,
\end{equation*}
where $R_1$ is defined in \eqref{def_R1}. In addition, if (A3) holds, then the solution to \eqref{weak_kirchhoff} is unique.
\end{theorem}

\begin{remark}
The condition $(m_0-L_f\lambda_1)>0$ guarantees that the elliptic diffusion dominates the linear part of the nonlinearity of $f$ in $u$. It ensures the positivity of the Galerkin operator on sufficiently large spheres, leading to the required a priori estimates. Consequently, it is essential for establishing the existence (and uniqueness) via the Brouwer fixed-point argument.
\end{remark}

\begin{remark}\label{rem_exit_Browder}
For the following particular type of nonlocal Kirchhoff problem, where $m_0=1$ and $L_f=0$ for $f(x,u)=f(x)$:
\begin{align}\label{particular_kirchhoff}
    -(1+\|\nabla u\|^2)\Delta u = f(x)\;\;\text{ in }\; \Omega,\;\; u=0\;\;\text{ in }\;\partial\Omega.
\end{align}
Assumptions (A1)--(A3) can be relaxed to prove the existence and uniqueness of the solution of \eqref{particular_kirchhoff} by the Browder-Minty theory; see \cite[Theorem~2.2 and Remark~2.3]{DondPani17_Kirchhoff} for a detailed proof.
\end{remark}

\section{Hybrid High-Order discretization}\label{sect:HHO_Dis}

\subsection{Discrete setting}
Let $\mathcal{H}\subset (0,\infty)$ be a countable set of meshsizes that has $0$ as its unique accumulation point. For $h\in\mathcal{H}$, we partition the domain $\Omega$ into a (polytopal) mesh $\mathcal{M}_h=(\cT_h,\cF_h)$; we refer \cite[Definition~1.4]{Piet_Jero_HHO_Book_20} for a detailed definition. The set of mesh elements $\cT_h$ consists of a finite collection of nonempty, disjoint polytopes $T$ with boundary $\partial T$ and diameter $h_T$ such that $\overline{\Omega}=\cup_{T\in\cT_h} \overline{T}$ and the meshsize satisfies $h=\max_{T\in\cT_h} h_T$. The set of mesh faces $\cF_h$ forms a partition of the mesh skeleton, that is, $\cup_{T\in\cT_h}\partial T=\cup_{F\in\cF_h}\bar{F}$, where $\cF_h=\cF_h^i\cup \cF_h^b$ with $\cF_h^i$ and $\cF_h^b$ representing the set of internal and boundary faces, respectively. 
For $F\in \cF_h$, we denote by $h_F$ the diameter of $F$. For each $T\in \cT_h$, we let $\mathcal{F}_T:=\{F\in \cF_h\, |\, F\subset\partial T\}$ be the set of faces belonging to $\partial T$. We denote $\bfn_T$ the unit outward normal to $T$ and set $\bfnTF:=\bfn_T|_F$ for all $F\in \cF_h$. We assume that the sequence of meshes $(\cM_h)_{h\in\mathcal{H}}$ is a regular mesh sequence in the sense of \cite[Assumption 1]{Jerome_Liam_21_hTscaling}, i.e., there exists a constant $\varrho>0$ such that for each $h\in\mathcal{H}$, each $T\in\cT_h$ and each $F\in\cF_h$ are connected by star-shaped sets with parameter $\varrho$ (see also \cite[Definition 1.41]{Piet_Jero_HHO_Book_20}).

% \noindent  Let $\mathbb{P}_d^l(T)$ be the polynomial space of degree at most $l$ on $T\in\cT_h$. There exist real numbers $\Ctr$ and $\Ctrc$ depending on $\varrho$ but independent of $h$ such that the following discrete and continuous trace inequalities hold for all $T\in\cT_h$ and $F\in\cF_T$ (see \cite[Lemma~1.46 and 1.49]{Piet_Ern_12_DG_book})
% \begin{align}
% 	\|v\|_F&\leq \Ctr h_F^{-1/2} \|v\|_T\quad\forall v\in\mathbb{P}_d^l(T),\label{dis_trace}\\ 
% 	\|v\|_{\partial T}&\leq C_{\text{tr},\text{c}}(h_T^{-1}\|v\|_T^2+ h_T\|\nabla v\|_T^2)^{1/2}\quad\forall v\in H^1(T).\label{cts_trace_ineq}
% \end{align}

Let $\mathbb{P}_d^l(T)$ be the polynomial space of degree at most $l$ on $T\in\cT_h$. Let $\pi_T^l:L^2(T)\to \mathbb{P}_d^l(T)$ denote the $L^2$-orthogonal projector. There exists a real number $C_{\text{app}}$ depending on $\varrho$ and $l$ but independent of $h$ such that the following approximation estimate holds (see \cite[Lemma~1.58 \& 1.59]{Piet_Ern_12_DG_book}): For all $ s\in \{1,\ldots,l+1\}$ and all $v\in H^s(T)$,
\begin{align}
	|v-\pi_T^l v|_{H^m(T)}+ h_T^{1/2}|v-\pi_T^l v|_{H^m(\partial T)}\leq C_{\text{app}} h_T^{s-m}|v|_{H^s(T)},\quad \forall m\in \{0, \ldots,s-1\},\label{proj_est}
\end{align}
where $|\bullet|_{H^m(\partial T)}$ denotes the facewise $H^m$-seminorm computed on $\partial T= \cup_{F\in \mathcal{F}_T}F$.

\subsection{Discrete spaces}
Let $k\geq 0$ denote a fixed polynomial degree. For a cell $T\in\cT_h$, let $\mathbb{P}_d^k(T)$ denote the space of polynomials of degree at most $k$ on $T$, and for a face $F\in\cF_h$, let $\mathbb{P}_{d-1}^k(F)$ denote the space of polynomial of degree at most $k$ on $F$. The local space of degrees of freedom (DOFs) on $T$ is defined as
\begin{align}\label{local_dofs}
	\underline{U}_T^k&:= \mathbb{P}_d^k(T)\times \left\{\underset{F\in\mathcal{F}_T}{\times} \mathbb{P}_{d-1}^k(F)\right\}.
\end{align}
The global space of DOFs is obtained by patching interface values in \eqref{local_dofs} as
\begin{align*}
	\underline{U}_h^k&:= \left\{\underset{T\in\mathcal{T}_h}{\times} \mathbb{P}_d^k(T)\right\}\times \left\{\underset{F\in\mathcal{F}_h}{\times} \mathbb{P}_{d-1}^k(F)\right\}.
\end{align*}
Imposing a homogeneous Dirichlet boundary condition on the above discrete space $\underline{U}_h^k$, we define
\begin{align*}
	\underline{U}_{h,0}^k&:= \left\{\underline{v}_h=\left((v_T)_{T\in\cT_h},(v_F)_{F\in\cF_h}\right)\in\underline{U}_h^k  \,|\, v_F\equiv 0\fl F\in\mathcal{F}_h^b\right\}.
\end{align*}
Let $\pi_F^k:L^2(F)\to \mathbb{P}_{d-1}^k(F)$ denote the $L^2$-orthogonal projector. Define a local interpolation operator $\underline{I}_T^k:H^1(T)\to \underline{U}_T^k$ such that for all $v\in H^1(T)$,
\begin{align}
	\underline{I}_T^kv:= (\pi_T^k v, (\pi_F^kv)_{F\in\mathcal{F}_T}).\label{defn_interpolant}
\end{align}
 The corresponding global interpolation operator $\underline{I}_h^k:H^1(\Omega)\to \underline{U}_h^k$ is given by
\begin{align*}
	\underline{I}_h^kv:= ((\pi_T^kv)_{T\in\mathcal{T}_h}, (\pi_F^kv)_{F\in\mathcal{F}_h})\fl v\in H^1(\Omega),
\end{align*}
and it maps $H_0^1(\Omega)$ onto $\underline{U}_{h,0}^k$.

\subsection{Local reconstructions}
For each $T\in\cT_h$, the local reconstruction operator $R_T^{k+1}:\underline{U}_T^k\to \mathbb{P}_d^{k+1}(T)$ is defined such that for $\underline{v}_T=(v_T,(v_F)_{F\in\cF_T})\in \underline{U}_T^k$,
\begin{subequations}\label{recons_oper}
	\begin{align}
		(\nabla R_T^{k+1}\underline{v}_T,\nabla w)_T&=(\nabla v_T,\nabla w)_T+ \sum_{F\in\mathcal{F}_T}(v_F-v_T,\nabla w{\cdot}\bfnTF)_F,\quad\forall w\in \mathbb{P}_d^{k+1}(T),\label{hho3}\\
		\left(R_T^{k+1}\underline{v}_T,1\right)_T&= \left(v_T,1\right)_T.\label{hho9}
	\end{align}
\end{subequations}
The corresponding global reconstruction operator $R_h^{k+1}:\underline{U}_h^k\to\mathbb{P}_d^{k+1}(\cT_h)$ is defined by $R_h^{k+1} \underline{v}_h|_T=R_T^{k+1} \underline{v}_T$.

The next lemma follows from \cite[Theorem~1.48]{Piet_Jero_HHO_Book_20} with the trace inequality and the approximation properties of an elliptic projector $\pi_T^{1,k+1}$ since $R_T^{k+1}\underline{I}_T^k v=\pi_T^{1,k+1}v$ for $v\in W^{1,1}(T)$.
\begin{lemma}[Approximation properties of $R_T^{k+1}\underline{I}_T^k$]\label{lem_apprx_recons}
 There exists a real number $C>0$, depending on $\varrho$ but independent of $h_T$ such that for all $v\in H^{s+1}(T)$ for some $s\in\{0,1,\ldots,k+1\}$, 
\begin{align}
&\|v-R_T^{k+1}\underline{I}_T^k v\|_T+ h_T^{1/2}\|v-R_T^{k+1}\underline{I}_T^k v\|_{\partial T}+ h_T\|\nabla(v-R_T^{k+1}\underline{I}_T^k v)\|_T \leq Ch_T^{s+1} |v|_{H^{s+1}(T)}.\label{eqn_recons_est}
\end{align}
For $s \in \{1,2,...,k+1\}$ and $v \in H^{s+1}(T)$, we also have the approximation property
\begin{align}\label{recon_approximation_1}
h_T^{1/2}\|\nabla(v-R_T^{k+1}\underline{I}_T^k v)\|_{\partial T}\leq Ch_T^{s} |v|_{H^{s+1}(T)}.
\end{align}
\end{lemma}

\section{Discrete HHO formulation}\label{sect:dis_HHO_form}
In this section, we design a hybrid high-order approximation for \eqref{cts_kirchhoff} and establish the existence of a unique discrete solution.
The discrete hybrid high-order formulation of \eqref{weak_kirchhoff} reads: Find $\underline{u}_h\in \underline{U}_{h,0}^k$ such that
\begin{equation}\label{dis_kirchhoff}
    M\left(\|\nabla R_h^{k+1}\underline{u}_h\|^2\right)a_h(\underline{u}_h,\underline{v}_h)=(f(u_h),v_h)\quad\forall \underline{v}_h\in\underline{U}_{h,0}^k,
\end{equation}
where 
\begin{equation}
a_h(\underline{u}_h,\underline{v}_h):=\sum_{T\in\cT_h}\left(a_T(\underline{u}_T,\underline{v}_T)+s_T(\underline{u}_T,\underline{v}_T)\right) \text{ with } a_T(\underline{u}_T,\underline{v}_T)=(\nabla R_T^{k+1}\underline{u}_T,\nabla R_T^{k+1}\underline{v}_T)_T
\end{equation}
and $$s_T(\underline{u}_T,\underline{v}_T):=\frac{1}{h_T}\sum_{F \in \cF_T} \left(\pi_{F}^{k}(u_F-u_T-(R_T^{k+1}\underline{u}_T-\pi_T^{k}R_T^{k+1}\underline{u}_T)),\pi_{F}^{k}(v_F-v_T-(R_T^{k+1}\underline{v}_T-\pi_T^{k}R_T^{k+1}\underline{v}_T))\right)_F.$$

An equivalent discrete formulation reads: Find $\underline{u}_h\in \underline{U}_{h,0}^k$ and $d\in\mR$ such that
\begin{align}
    \|\nabla R_h^{k+1}\underline{u}_h\|^2-d&=0\label{equivalent_dis1}\\
    M(d)a_h(\underline{u}_h,\underline{v}_h)-(f(u_h),v_h)&=0 \quad\forall \underline{v}_h\in\underline{U}_{h,0}^k.\label{equivalent_dis2}
\end{align}
We define the discrete norms on $\underline{U}_{h,0}^k$ as: for $\underline{v}_h\in \underline{U}_{h,0}^k$
\begin{equation*}
    \|\underline{v}_h\|_{a,h}^2:=a_h(\underline{v}_h,\underline{v}_h)\text{ and }\|\underline{v}_h\|_{1,h}^2:=\sum_{T\in\cT_h}\big{(}\|\nabla v\|_T^2+\sum_{F\in\cF_T}h_F^{-1}\|v_F-v_T\|_F^2\big{)}.
\end{equation*}
Note that the above energy norm $\|\cdot\|_{a,h}$ is equivalent to the standard HHO norm $\|\cdot\|_{1,h}$; see \cite[Lemma~4]{Piet_Ern_Lem_14_arb_local}.
The next theorem on the equivalence of formulations \eqref{dis_kirchhoff} and \eqref{equivalent_dis1}-\eqref{equivalent_dis2} follows immediately.

\begin{theorem}
    If $(\underline{u}_h,d)\in \underline{U}_{h,0}^k\times \mR$ is a solution of \eqref{equivalent_dis1}-\eqref{equivalent_dis2}, then $\underline{u}_h$ is a solution to \eqref{dis_kirchhoff}. Conversely, if $\underline{u}_h \in \underline{U}_{h,0}^k$ is a solution of \eqref{dis_kirchhoff}, then $(\underline{u}_h,d):=(\underline{u}_h,\|\nabla R_h^{k+1}\underline{u}_h\|^2)$ is a solution to \eqref{equivalent_dis1}-\eqref{equivalent_dis2}.
\end{theorem}

The next result follows from the Brouwer's fixed point theorem on Euclidean space; see \cite[Theorem~5.2.5]{Kesavan2008}.
\begin{proposition}\label{prop_fixed_point}
Let $H$ be a finite dimensional Hilbert space with inner product $(\cdot,\cdot)$ and norm $|\cdot|$. Let $\mathcal{S}:H\to H$ be a continuous map with the following property: there exists $\rho>0$ such that
\begin{equation*}
     (\mathcal{S}(v),v)\geq 0\quad\forall v\in H\text{ with } |v|=\rho.
\end{equation*}
Then, there exists a $w\in H$ such that $\mathcal{S}(w)=0$ with $|w|\leq \rho$. 
\end{proposition}

% For $\underline{u}_h,\underline{v}_h\in \underline{U}_{h,0}^k$, define
% \begin{align}
%     B(\underline{u}_h;\underline{u}_h,\underline{v}_h)&:=M\left(\|\nabla R_h^{k+1}\underline{u}_h\|^2\right)a_h(\underline{u}_h,\underline{v}_h),\label{def_B}\\
%     L(\underline{u}_h;\underline{v}_h)&:=(f(u_h),v_h).\label{def_L}
% \end{align}

% To obtain an estimate for a discrete approximate eigenvalue, we consider the following spectral HHO approximation problem of \cite{Cicu_Ern_19_Eigen_HHO}: find $\underline{u}_h\in\underline{U}_{h,0}^k$ such that
% \begin{equation}\label{dis_eigen}
%     a_h(\underline{u}_h,\underline{v}_h)=\lambda_h (u_h,v_h)\quad\forall \underline{v}_h\in\underline{U}_{h,0}^k.
% \end{equation}
% From \cite[Theorem~4.4]{Cicu_Ern_19_Eigen_HHO}, we obtain an estimate for the discrete principal eigenvalue $\lambda_{1h}^{-1}$ of \eqref{dis_eigen} as 
% \begin{equation}\label{eigen_est}
%     |\lambda_1^{-1}-\lambda_{1h}^{-1}|=\mathcal{O}(h^{2t}), 
% \end{equation}
% for some $t>\frac{1}{2}$ that depends on the regularity of the solution. The above equations \eqref{dis_eigen} and \eqref{eigen_est} lead to 
% \begin{equation}\label{dis_eigen_bdd}
%     \|u_h\|^2\leq \lambda_{1h}^{-1} a_h(\underline{u}_h,\underline{u}_h)\leq\left(\lambda_1^{-1}+\mathcal{O}(h^{2t})\right)\|\underline{u}_h\|_{a,h}^2.
% \end{equation}

Let $\lambda_{1,h}$ be the smallest (discrete) eigenvalue of the discrete eigenvalue problem; see \cite{Cicu_Ern_19_Eigen_HHO}: Find
$\underline{w}_h\in\underline{U}_{h,0}^k$ such that
\begin{equation}\label{dis_eigen}
    a_h(\underline{w}_h,\underline{v}_h)=\lambda_h (w_h,v_h)\quad\forall \underline{v}_h\in\underline{U}_{h,0}^k.
\end{equation}
For nonconforming (or discontinuous Galerkin) finite element methods, inequality $\lambda_{1,h}\geq \lambda_1$ cannot, in general, be guaranteed. Consequently, an $L^2$-estimate of discrete HHO functions in terms of $\lambda_1^{-1/2}\|\cdot\|_{a,h}$ (as in Gudi \cite{Gudi12_Kirchhof_apriori} for the conforming FEM for the nonlocal problem of Kirchhoff type) is not immediately available.
However, using the finite-dimensional Rayleigh principle and the spectral convergence $\lambda_{1,h}=\lambda_1 + \mathcal{O}(h^{2t})$ with $t> 1/2$ of \cite[Theorem~4.4]{Cicu_Ern_19_Eigen_HHO}, the following result holds:

\begin{lemma}\label{lmm_l2est} There exist positive constants $C_R$ and $h_0$ such that, for all $0<h\leq h_0$ and for all $\underline{v}_h\in\underline{U}_{h,0}^k$,
\begin{equation}\label{dis_eigen_bdd}
\|v_h\|^2\leq \lambda_{1,h}^{-1} a_h(\underline{v}_h,\underline{v}_h)\leq\left(\lambda_1^{-1}+C_Rh^{2t})\right)\|\underline{v}_h\|_{a,h}^2,
\end{equation}
where $t\in (1/2,k+1]$ is the elliptic regularity index.
\end{lemma}
% \begin{proof}
% By the finite–dimensional Rayleigh principle, the smallest discrete eigenvalue satisfies
% \begin{equation*}
%     \lambda_{1,h}=\min_{\underline{v}_h\in\underline{U}_{h,0}^k\setminus\{0\}}\frac{a_h(\underline{v}_h,\underline{v}_h)}{(v_h,v_h)}.
% \end{equation*}
% This implies $\|v_h\|^2\leq \lambda_{1,h}^{-1} a_h(\underline{v}_h,\underline{v}_h)$ for all $\underline{v}_h\in\underline{U}_{h,0}^k$. The spectral convergence $\lambda_{1,h}=\lambda_1 + \mathcal{O}(h^{2t})$ of \cite[Theorem~4.4]{Cicu_Ern_19_Eigen_HHO} implies $\lambda_{1,h}^{-1}=\lambda_1^{-1} + \mathcal{O}(h^{2t})$. Therefore, there exist constants $C_R$ and $h_0>0$ such that $\lambda_{1,h}^{-1}\leq \lambda_1^{-1} + C_R h^{2t}$. This completes the proof.
% \end{proof}

% \begin{lemma}
% Let $\underline{u}_h,\underline{v}_h\in\underline{U}_{h,0}^k$, then
% \begin{equation*}
% \left|M(l_h(\underline{u}_h))-M(l_h(\underline{v}_h))\right|\leq 
% \begin{cases}
%     L_M\|\underline{u}_h-\underline{v}_h\|_{a,h}\left(\|\nabla R_h^{k+1}\underline{u}_h\|+\|\nabla R_h^{k+1}\underline{v}_h\|\right)\text{ if } l_h(\underline{u}_h) = \|\nabla R_h^{k+1}\underline{u}_h\|^2\\
%     L_M|\Omega|^{1/2}\lambda_{1,h}^{-1/2}\|\underline{u}_h-\underline{v}_h\|_{a,h} \qquad\qquad\qquad\quad\text{ if } l_h(\underline{u}_h) = \int_{\Omega} u_h\dx.
% \end{cases}
% \end{equation*}
% \end{lemma}
Now we prove the existence and uniqueness of the solution to \eqref{dis_kirchhoff}.
\begin{theorem}
Suppose that (A1)--(A2) hold. Then for sufficiently small $h$, the discrete HHO problem~\eqref{dis_kirchhoff} has a solution $\underline{u}_h\in\underline{U}_{h,0}^k$. Moreover, any solution $\underline{u}_h$ of \eqref{dis_kirchhoff} satisfies $\|\underline{u}_h\|_{a,h}\leq R_{1,h}$, where 
$$R_{1,h}:=\frac{\|f(\cdot,0)\|}{\lambda_{1,h}^{1/2}(m_0-L_f\lambda_{1,h}^{-1})}.$$ In addition, if (A3) holds, the problem \eqref{dis_kirchhoff} has a unique solution. Consequently, the equivalent problem~\eqref{equivalent_dis1}-\eqref{equivalent_dis2} has a unique solution.
\end{theorem}
\begin{proof}
For a fixed $\underline{w}_h\in \underline{U}_{h,0}^k$, define a linear map $Q_{w}:\underline{U}_{h,0}^k\to \mR$ by 
\begin{equation}\label{defn_Qw}
    Q_{w}(\underline{v}_h):=M\left(\|\nabla R_h^{k+1}\underline{w}_h\|^2\right)a_h(\underline{w}_h,\underline{v}_h)-(f(w_h),v_h).
\end{equation}
The Cauchy-Schwarz inequality leads to
\begin{equation}
    |M\left(\|\nabla R_h^{k+1}\underline{w}_h\|^2\right)a_h(\underline{w}_h,\underline{v}_h)|\leq M\left(\|\nabla R_h^{k+1}\underline{w}_h\|^2\right)\|\underline{w}_h\|_{a,h}\|\underline{v}_h\|_{a,h}.\label{bd_B}
\end{equation}
The Cauchy-Schwarz inequality and Lipschitz continuity of $f$ lead to
\begin{align}
    |(f(w_h),v_h)|&\leq\int_{\Omega}|f(w_h)||v_h|\dx\leq \int_{\Omega}(L_f|w_h|+|f(x,0)|)|v_h|\dx\notag\\
    &\leq \|L_f|w_h|+|f(x,0)|\|\|v_h\|.\label{bd_L}
\end{align}
The above two estimates \eqref{bd_B}-\eqref{bd_L} lead to a boundedness result of $Q_{w}$ as
\begin{equation}
    |Q_{w}(\underline{v}_h)|\leq \left( M\left(\|\nabla R_h^{k+1}\underline{w}_h\|^2\right)\|\underline{w}_h\|_{a,h}+\lambda_{1,h}^{-1/2}\|L_f|w_h|+|f(x,0)|\|\right)\|\underline{v}_h\|_{a,h}.\label{Q_bdd}
\end{equation}
Therefore, $Q_{w}$ is a continuous linear map on $\underline{U}_{h,0}^k$ equipped with the norm $\|\cdot\|_{a,h}$. By Riesz representation theorem, there exists a unique $\mathcal{S}(\underline{w}_h)\in \underline{U}_{h,0}^k$ such that
\begin{equation}\label{defn_S}
    a_h(\mathcal{S}(\underline{w}_h),\underline{v}_h)=Q_{w}(\underline{v}_h)\quad\forall\underline{v}_h\in \underline{U}_{h,0}^k.
\end{equation}
Since $M$ is continuous and $f$ is Lipschitz continuous, the nonlinear map $\mathcal{S}:\underline{U}_{h,0}^k\to \underline{U}_{h,0}^k$ is continuous; see appendix for a detailed proof.
% Now we show that there exists a $\underline{u}_h\in\underline{U}_{h,0}^k$ such that
% \begin{equation}
%     a_h(R(\underline{u}_h),\underline{v}_h)=0\quad\forall\underline{v}_h\in \underline{U}_{h,0}^k.
% \end{equation}
% This will prove that problem~\eqref{dis_kirchhoff} has a solution in $\underline{U}_{h,0}^k$. 
% {\color{red}From the boundedness of $Q_w$ in \eqref{Q_bdd}, we have 
% \begin{align}
%     \|\mathcal{S}(\underline{w}_h)\|_{a,h}&=\sup_{\|\underline{v}_h\|_{a,h}\leq 1}a_h(\mathcal{S}(\underline{w}_h),\underline{v}_h)=\sup_{\|\underline{v}_h\|_{a,h}\leq 1}Q_{w}(\underline{v}_h)\notag\\
%     &\leq  M\left(\|\nabla R_h^{k+1}\underline{w}_h\|^2\right)\|\underline{w}_h\|_{a,h}+\lambda_{1,h}^{-1/2}L_f\|w_h\|+\|f(x,0)\|.
% \end{align}
% Since $M(\cdot)$ is continuous, the nonlinear map $\mathcal{S}:\underline{U}_{h,0}^k\to \underline{U}_{h,0}^k$ is continuous.} 
Now we prove the positivity of $\mathcal{S}$ as follows. Using the definition of $Q_w$, Assumption (A2) and \eqref{dis_eigen_bdd}, we have
\begin{align}
    a_h(\mathcal{S}(\underline{w}_h),\underline{w}_h)&=Q_{w}(\underline{w}_h)\notag\\
    &=M\left(\|\nabla R_h^{k+1}\underline{w}_h\|^2\right)a_h(\underline{w}_h,\underline{w}_h)-(f(w_h),w_h)\notag\\
    &\geq m_0\|\underline{w}_h\|_{a,h}^2-L_f\|w_h\|^2-\|f(\cdot,0)\|\|w_h\|\notag\\
    &\geq (m_0-L_f\lambda_{1,h}^{-1})\|\underline{w}_h\|_{a,h}^2-\lambda_{1,h}^{-1/2}\|f(\cdot,0)\|\|\underline{w}_h\|_{a,h}.
\end{align}
Therefore, $a_h(\mathcal{S}(\underline{w}_h),\underline{w}_h)\geq 0$ for all $\underline{w}_h\in \underline{U}_{h,0}^k$ with $\|\underline{w}_{h}\|_{a,h}=\rho$, where $\rho$ is any positive number such that 
\begin{equation}
    \rho\geq\frac{\lambda_{1,h}^{-1/2}\|f(\cdot,0)\|}{(m_0-L_f\lambda_{1,h}^{-1})}=R_{1,h}.
\end{equation}
This with Proposition~\ref{prop_fixed_point} completes the proof of the existence of a solution $\underline{u}_h\in \underline{U}_{h,0}^k$ of \eqref{dis_kirchhoff} with $\|\underline{u}_h\|_{a,h}\leq R_{1,h}$.

To prove the uniqueness, we make use of Assumption (A3). Suppose that there are two solutions $\underline{u}_h^1$ and $\underline{u}_h^2$ of \eqref{dis_kirchhoff} such that $\|\underline{u}_h^1\|_{a,h}\leq R_{1,h}$ and $\|\underline{u}_h^2\|_{a,h}\leq R_{1,h}$. Let $\underline{\psi}_h:=\underline{u}_h^1-\underline{u}_h^2$. Then, (A2) and simple arrangement of terms lead to
\begin{align}
    m_0\|\underline{\psi}_h\|_{a,h}^2&\leq M\left(\|\nabla R_h^{k+1}\underline{u}_h^1\|^2\right)\|\underline{\psi}_h\|_{a,h}^2\notag\\
    &=M\left(\|\nabla R_h^{k+1}\underline{u}_h^1\|^2\right)a_h(\underline{u}_h^1,\underline{\psi}_h)-M\left(\|\nabla R_h^{k+1}\underline{u}_h^1\|^2\right)a_h(\underline{u}_h^2,\underline{\psi}_h)\notag\\
    &=\left(M\left(\|\nabla R_h^{k+1}\underline{u}_h^1\|^2\right)-M\left(\|\nabla R_h^{k+1}\underline{u}_h^2\|^2\right)\right)a_h(\underline{u}_h^2,\underline{\psi}_h)+\left(f(u_h^1)-f(u_h^2),\psi_h\right).\label{est_coer}
\end{align}
Since $M$ is Lipschitz continuous and $\|\nabla R_h^{k+1}\underline{\psi}_h\|\leq \|\underline{\psi}_h\|_{a,h}$, we have the following estimate:
\begin{align}
    &\left|M\left(\|\nabla R_h^{k+1}\underline{u}_h^1\|^2\right)-M\left(\|\nabla R_h^{k+1}\underline{u}_h^2\|^2\right)\right|\leq L_M\left|\|\nabla R_h^{k+1}\underline{u}_h^1\|^2-\|\nabla R_h^{k+1}\underline{u}_h^2\|^2\right|\notag\\
    &\leq L_M\|\nabla R_h^{k+1}(\underline{u}_h^1-\underline{u}_h^2)\|(\|\nabla R_h^{k+1}\underline{u}_h^1\|+\|\nabla R_h^{k+1}\underline{u}_h^2\|)\notag\\
    &\leq 2L_M R_{1,h}\|\nabla R_h^{k+1}\underline{\psi}_h\|\leq 2L_M R_{1,h}\| \underline{\psi}_h\|_{a,h}.\label{est_M}
\end{align}
Since $f$ is also Lipschitz continuous, we have
\begin{align}
  \left|\left(f(u_h^1)-f(u_h^2),\psi_h\right)\right|\leq L_f\|u_h^1-u_h^2\|\|\psi_h\|\leq L_f \|\psi_h\|^2\leq L_f\lambda_{1,h}^{-1}\|\underline{\psi}_h\|_{a,h}^2.\label{est_f}
\end{align}
Combining the estimates \eqref{est_M} and \eqref{est_f} in \eqref{est_coer}, we obtain
\begin{equation}
    \left(m_0-2L_M R_{1,h}^2-L_f\lambda_{1,h}^{-1}\right)\|\underline{\psi}_h\|_{a,h}^2\leq 0.
\end{equation}
We note that $R_{1,h}\to R_{1}$ and $\lambda_{1,h}\to \lambda_1$ as $ h\to 0$. Since $m_0-2L_M R_1^2-L_f\lambda_1^{-1}>0$, for a sufficiently small mesh parameter $h$, we have $m_0-2L_M R_{1,h}^2-L_f\lambda_{1,h}^{-1}>0$. Therefore, we obtain $\|\underline{\psi}_h\|_{a,h}=0$, that is, $\underline{u}_h^1=\underline{u}_h^2$. This completes the proof of the uniqueness of the solution of \eqref{dis_kirchhoff}.
\end{proof}

\begin{remark}
For the particular type of Kirchhoff problem~\eqref{particular_kirchhoff}, we can prove the existence of a unique discrete HHO solution without assumptions (A1)--(A3). However, we need to consider the following (stabilized coefficient) discrete reformulation: Find $\underline{u}_h\in \underline{U}_{h,0}^k$ such that
\begin{equation}\label{dis_kirchhoff_specific}
    \left(1+a_h(\underline{u}_h,\underline{u}_h)\right)a_h(\underline{u}_h,\underline{v}_h)=(f,v_h)\quad\forall \underline{v}_h\in\underline{U}_{h,0}^k.
\end{equation}
Define $\mathcal{A}: \underline{U}_{h,0}^k\to \underline{U}_{h,0}^k$ by
$$a_h(\mathcal{A}(\underline{u}_h),\underline{v}_h):=\left(1+a_h(\underline{u}_h,\underline{u}_h)\right)a_h(\underline{u}_h,\underline{v}_h)\quad\forall \underline{v}_h\in\underline{U}_{h,0}^k.$$
The operator satisfies the monotonicity property (see also \cite[Theorem~2.2]{DondPani17_Kirchhoff}): for all $\underline{v}_h,\underline{w}_h\in\underline{U}_{h,0}^k$
\begin{align}
a_h(\mathcal{A}(\underline{v}_h)-\mathcal{A}(\underline{w}_h),\underline{v}_h-\underline{w}_h)&=a_h(\underline{v}_h-\underline{w}_h,\underline{v}_h-\underline{w}_h)\notag\\
&\quad+\left(a_h(\underline{v}_h,\underline{v}_h)a_h(\underline{v}_h,\underline{v}_h-\underline{w}_h)-a_h(\underline{w}_h,\underline{w}_h)a_h(\underline{w}_h,\underline{v}_h-\underline{w}_h)\right)\notag\\
&\geq a_h(\underline{v}_h-\underline{w}_h,\underline{v}_h-\underline{w}_h) + \frac{1}{2}\left(\|\underline{v}_h\|_{a,h}^2-\|\underline{w}_h\|_{a,h}^2\right)^2\notag\\
&\geq a_h(\underline{v}_h-\underline{w}_h,\underline{v}_h-\underline{w}_h).
\end{align}
The existence and uniqueness of the solution of \eqref{dis_kirchhoff_specific} follow from the above monotonicity property and an application of Browder-Minty theorey; see \cite[Remark~2.3]{DondPani17_Kirchhoff}.
\end{remark}

\section{Error estimate}\label{sect:err_est}
In this section, we derive an error estimate of the discrete HHO solution.
\begin{theorem}\label{thm:Kirchhoof_err_est}
    Suppose that (A1)--(A3) hold. Assume that the solution of \eqref{weak_kirchhoff} has the regularity $u\in H^1_0(\Omega)\cap H^{r+2}(\cT_h)$ for some $r\in\{0,\ldots,k\}$. Then, for sufficiently small $h$, the discrete HHO solution $\underline{u}_h\in \underline{U}_{h,0}^k$ of \eqref{dis_kirchhoff} satisfies the error estimate
    \begin{equation}\label{eqn_err_est}
        \|\underline{I}_h^ku-\underline{u}_h\|_{a,h}\leq Ch^{r+1}\|u\|_{H^{r+2}(\cT_h)}.
    \end{equation}
\end{theorem}
\begin{proof}
Let $\underline{\psi}_h=\underline{I}_h^ku-\underline{u}_h$. Then, using \eqref{cts_kirchhoff} and \eqref{dis_kirchhoff}, we express
\begin{align}
M(\|\nabla u\|^2)a_h(\underline{\psi}_h,\underline{\psi}_h)&=M(\|\nabla u\|^2)a_h(\underline{I}_h^k u,\underline{\psi}_h)-M(\|\nabla u\|^2)a_h(\underline{u}_h,\underline{\psi}_h)\notag\\
&=M(\|\nabla u\|^2)\left(a_h(\underline{I}_h^k u,\underline{\psi}_h)+(\Delta u,\psi_h)\right)+(f(u),\psi_h)-M(\|\nabla u\|^2)a_h(\underline{u}_h,\underline{\psi}_h)\notag\\
&=M(\|\nabla u\|^2)\left(a_h(\underline{I}_h^k u,\underline{\psi}_h)+(\Delta u,\psi_h)\right)+\left(f(u)-f(\pi_h^k u),\psi_h\right)\notag\\
&\quad+\left(f(\pi_h^k u)-f(u_h),\psi_h\right)+\left(M(\|\nabla R_h^{k+1}\underline{u}_h\|^2)-M(\|\nabla u\|^2)\right)a_h(\underline{u}_h,\underline{\psi}_h).\label{est_err}
\end{align}
The first term is related to the general consistency error,  $\mathcal{E}_h(w;\underline{\psi}_h):=a_h(\underline{I}_h^k w,\underline{\psi}_h)+(\Delta w,\psi_h)$ of \cite[Lemma~2.18]{Piet_Jero_HHO_Book_20} with convergence 
\begin{equation}\label{conv_consistency}
    |\mathcal{E}_h(w;\cdot)|\leq C h^{r+1}|w|_{H^{r+2}(\cT_h)},
\end{equation}
for some positive constant $C$ independent of $h$.
Therefore, the first term is written as
\begin{equation}
 M(\|\nabla u\|^2)\left(a_h(\underline{I}_h^k u,\underline{\psi}_h)+(\Delta u,\psi_h)\right)\leq M(\|\nabla u\|^2)  \mathcal{E}_h(u;\underline{\psi}_h).\label{est_consist}   
\end{equation}
The second term of \eqref{est_err} is estimated using the Lipschitz continuity of $f$ as
\begin{equation}
 \left|\left(f(u)-f(\pi_h^k u),\psi_h\right)\right|\leq L_f\|u-\pi_h^k u\|\|\psi_h\|\leq L_f\lambda_{1,h}^{-1/2}\|u-\pi_h^k u\|\|\underline{\psi}_h\|_{a,h}.
\end{equation}
The third term of \eqref{est_err} has the estimate
\begin{equation}
    \left|\left(f(\pi_h^k u)-f(u_h),\psi_h\right)\right|\leq L_f \|\pi_h^k u-u_h\|\|\psi_h\|\leq L_f\|\psi_h\|^2\leq L_f\lambda_{1,h}^{-1}\|\underline{\psi}_h\|_{a,h}^2.
\end{equation}
The last term of \eqref{est_err} is estimated using the Lipschitz continuity of $M$ with $\|\nabla R_h^{k+1}\underline{u}_h\|\leq \|\underline{u}_h\|_{a,h}\leq R_{1,h}$ and $\|\nabla R_h^{k+1}\underline{I}_h^k u\|\leq \|\nabla u\|\leq R_1$ as
\begin{align}
&\left(M(\|\nabla R_h^{k+1}\underline{u}_h\|^2)-M(\|\nabla u\|^2)\right)a_h(\underline{u}_h,\underline{\psi}_h)\notag\\
&\leq\left[\left(M(\|\nabla R_h^{k+1}\underline{u}_h\|^2)-M(\|\nabla R_h^{k+1}\underline{I}_h^k u\|^2)\right)+\left(M(\|\nabla R_h^{k+1}\underline{I}_h^k u\|^2)-M(\|\nabla u\|^2)\right)\right]a_h(\underline{u}_h,\underline{\psi}_h)\notag\\
&\leq L_M\Big[\left|\|\nabla R_h^{k+1}\underline{u}_h\|^2-\|\nabla R_h^{k+1}\underline{I}_h^k u\|^2\right|+\left|\|\nabla R_h^{k+1}\underline{I}_h^k u\|^2-\|\nabla u\|^2\right|\Big]|a_h(\underline{u}_h,\underline{\psi}_h)|\notag\\
&\leq L_M\Big[\|\nabla R_h^{k+1}(\underline{u}_h-\underline{I}_h^k u)\|(\|\nabla R_h^{k+1}\underline{u}_h\|+\|\nabla R_h^{k+1}\underline{I}_h^k u\|)\Big]|a_h(\underline{u}_h,\underline{\psi}_h)|\notag\\
&\quad+L_M\Big[\|\nabla R_h^{k+1}\underline{I}_h^k u-\nabla u\|(\|\nabla R_h^{k+1}\underline{I}_h^k u\|+\|\nabla u\|)\Big]|a_h(\underline{u}_h,\underline{\psi}_h)|\notag\\
&\leq L_M (R_{1,h}+R_1)R_{1,h}\|\underline{\psi}_h\|_{a,h}^2+2L_M R_1 R_{1,h} \|\nabla R_h^{k+1}\underline{I}_h^k u-\nabla u\|\|\underline{\psi}_h\|_{a,h}.\label{est_M_err}
\end{align}
Combining the above four estimates \eqref{est_consist}-\eqref{est_M_err} in \eqref{est_err}, we obtain
\begin{align}
&\left(M(\|\nabla u\|^2)-L_f\lambda_{1,h}^{-1}-L_M (R_{1,h}+R_1)R_{1,h}\right)\|\underline{I}_h^ku-\underline{u}_h\|_{a,h}\notag\\
&\leq M(\|\nabla u\|^2)  \mathcal{E}_h(u;\cdot)+L_f\lambda_{1,h}^{-1/2}\|u-\pi_h^k u\|+2L_M R_1 R_{1,h}\|\nabla( R_h^{k+1}\underline{I}_h^k u- u)\|.
\end{align}
We note that $R_{1,h}\to R_{1}$ and $\lambda_{1,h}\to \lambda_1$ as $ h\to 0$. Since $M(\|\nabla u\|^2)\geq m_0$ and $(m_0-L_f\lambda_1^{-1}-2L_M R_1^2)>0$ (from assumption (A3)), for a sufficiently small mesh parameter $h$, we have 
\begin{equation*}
    \left(m_0-L_f\lambda_{1,h}^{-1}-L_M (R_{1,h}+R_1)R_{1,h}\right)>0.
\end{equation*}
Using \eqref{conv_consistency}, \eqref{proj_est} and \eqref{eqn_recons_est}, we obtain the required error estimate
\begin{equation*}
 \|\underline{I}_h^ku-\underline{u}_h\|_{a,h}\leq Ch^{r+1}\|u\|_{H^{r+2}(\cT_h)},   
\end{equation*}
for some constant $C$ independent of the mesh parameter $h$.
\end{proof}

% \begin{remark}
% We can derive the following error estimate by adding and subtracting $\nabla R_h^{k+1} \underline{I}_h^k u$ as follows:
% \begin{align*}
% \|\nabla(u -R_h^{k+1} \underline{u}_h)\|&\leq \|\nabla(u -R_h^{k+1} \underline{I}_h^k u)\|+\|\nabla R_h^{k+1} (\underline{I}_h^k u - \underline{u}_h)\|\\
%         &\leq \|\nabla(u -R_h^{k+1} \underline{I}_h^k u)\|+\|\underline{I}_h^k u - \underline{u}_h\|_{a,h}.
% \end{align*}
% Applying Lemma 3.1 on the interpolation error and Theorem 5.1. on the discrete energy error estimate, we obtain the error estimate: $\|\nabla(u -R_h^{k+1} \underline{u}_h)\|\leq Ch^{r+1}.$
% \end{remark}
\begin{remark}
Although Theorem~\ref{thm:Kirchhoof_err_est} is stated in terms of the discrete energy norm, it immediately yields an error estimate for the reconstructed potential in the $H^1$-seminorm. Indeed, by adding and subtracting $R_h^{k+1}\underline{I}_h^k u$, we obtain
\begin{align*}
\|\nabla(u-R_h^{k+1}\underline{u}_h)\|
&\le
\|\nabla(u-R_h^{k+1}\underline{I}_h^k u)\|
+\|\nabla R_h^{k+1}(\underline{I}_h^k u-\underline{u}_h)\| \\
&\lesssim
\|\nabla(u-R_h^{k+1}\underline{I}_h^k u)\|
+\|\underline{I}_h^k u-\underline{u}_h\|_{a,h},
\end{align*}
where the second inequality follows from the stability of the potential reconstruction operator. The first term is bounded by the interpolation estimate of Lemma~\ref{lem_apprx_recons}, while the second term is estimated by Theorem~\ref{thm:Kirchhoof_err_est}. Consequently,
\[
\|\nabla(u-R_h^{k+1}\underline{u}_h)\|
\leq Ch^{r+1}\|u\|_{H^{r+2}(\cT_h)},
\]
where $r$ is as defined in Theorem~\ref{thm:Kirchhoof_err_est}.
\end{remark}
\begin{remark}
The error estimate \eqref{eqn_err_est} can also be derived using similar arguments for the particular Kirchhoff type problem~\eqref{particular_kirchhoff} and the corresponding discrete solution \eqref{dis_kirchhoff_specific} without assumptions (A1)--(A3).
\end{remark}

\section{Numerical Experiments}\label{sect:numer_exp}
In this section, we perform numerical tests to compute the discrete solution of \eqref{weak_kirchhoff} and illustrate the convergence rate. To solve the nonlocal nonlinear discrete problem, we use the Newton-Raphson method; see \cite{Gudi12_Kirchhof_apriori} for the implementation of a conforming finite element method. We describe the steps for the implementation of the discrete HHO solution \eqref{dis_kirchhoff} using the solution of the equivalent problem~\eqref{equivalent_dis1}-\eqref{equivalent_dis2}.  Let $N$ be the dimension of the HHO space $\underline{U}_{h,0}^k$ and $\{\underline{\varphi_i}\}_{i=1}^{N}$ be a basis of $\underline{U}_{h,0}^k$. Then, any solution $\underline{u}_h\in \underline{U}_{h,0}^k$ can be written as 
\begin{equation*}
    \underline{u}_h=\sum_{i=1}^{N}\alpha_i\underline{\varphi_i}.
\end{equation*}
Set $\vec{\alpha}=[\alpha_1,\alpha_2,\ldots,\alpha_N]^T$. Define a nonlinear map $F:\mR^{N}\times\mR\to\mR^{N+1}$ with $F=[F_1,F_2,\ldots,F_{N+1}]^T$ such that
\begin{align}
    F_j(\vec{\alpha},d)&=M(d)\sum_{i=1}^N\alpha_i a_h(\underline{\varphi_i},\underline{\varphi_j})-(f(u_h),\varphi_j), \text{ where } 1\leq j\leq N \text{ and } u_h = \sum_{i=1}^{N}\alpha_i\varphi_i,\\
    F_{N+1}(\vec{\alpha},d)&= \vec{\alpha}^T R\vec{\alpha}-d, \text{ where } R=(r_{ij})_{N\times N} \text{ with } r_{ij}=(\nabla R_h^{k+1}\underline{\varphi}_i,\nabla R_h^{k+1}\underline{\varphi}_j).
\end{align}
Now we compute the Jacobian matrix of the form
\begin{equation*}
   J(\vec{\alpha},d)= \begin{bmatrix}
        A & b\\
        c & \delta
    \end{bmatrix}
\end{equation*}
as follows:
\begin{align*}
    A_{ji} &= \frac{\partial F_{j}}{\partial \alpha_i}(\vec{\alpha},d) =M(d)a_h(\underline{\varphi_i},\underline{\varphi_j})-(f'(u_h)\varphi_i,\varphi_j), \quad 1\leq i,j\leq N,\\
    b_{j1}& =  \frac{\partial F_{j}}{\partial d}(\vec{\alpha},d) =M'(d)\sum_{i=1}^N\alpha_i a_h(\underline{\varphi_i},\underline{\varphi_j}), \quad 1\leq j\leq N,\\
    c_{1i}&=\frac{\partial F_{N+1}}{\partial \alpha_{i}}(\vec{\alpha},d) = 2(\nabla R_h^{k+1}\underline{\varphi}_i,\nabla R_h^{k+1}\underline{u}_h), \quad 1\leq i\leq N,\\
    \delta &=  \frac{\partial F_{N+1}}{\partial d}(\vec{\alpha},d) = -1.
\end{align*}
We observe that $A=A_{N\times N}$ is sparse and $b=b_{N\times 1}$ and $c=c_{1\times N}$ are full matrices. The matrix $A$ has the same sparsity structure as the Jacobian corresponding to the discrete HHO equation:
\begin{equation*}
    a_h(\underline{w}_h,\underline{v}_h)=(f(w_h),v_h)
\end{equation*}
of the continuous semilinear equation:
\begin{equation*}
(\nabla w,\nabla v)=(f(w),v).
\end{equation*}

For an initial guess $\beta^{(0)}=(\vec{\alpha}^{(0)},d^{(0)})$, sufficiently close to the exact solution $\beta:=(\vec{\alpha},d)$, Newton's iterations are given by 
\begin{align}\label{Newton_iter}
    J^{(n)}(\beta^{(n)})\left(\beta^{(n+1)}-\beta^{(n}\right)=-F^{(n)}(\beta^{(n)})\quad\text{with}\quad \beta^{(n)}:=(\vec{\alpha}^{(n)},d^{(n)}),
\end{align}
where $J^{(n)}$ and $F^{(n)}$ are computed at the $n$-th iteration using the $n$-th iteration approximate solution $\vec{\alpha}^{(n)}$ and $d^{(n)}$. We solve \eqref{Newton_iter} using the Sherman--Morrison-Woodbury formula (see \cite{Gudi12_Kirchhof_apriori} for details) and the static condensation of \cite{Piet_Ern_Lem_14_arb_local,Cicu_Piet_Ern_18_Comp_HHO}.
We compute the empirical rate of convergence for the relative $H^1$-error $e_h:=\|\underline{I}_h^ku-\underline{u}_h\|_{a,h}/\|\underline{I}_h^k u\|_{a,h}$ as follows:
\begin{align*}
\texttt{rate}(\ell):=\log \big(e_{h_{\ell}}/e_{h_{\ell-1}} \big)/\log \big(h_\ell/h_{\ell-1} \big)\text{ for } \ell=1,2,3,\ldots,
\end{align*}
where $e_{h_{\ell}}$ and $e_{h_{\ell-1}}$ are the errors associated with the two consecutive meshsizes $h_\ell$ and $h_{\ell-1}$, respectively.

\begin{figure}
	\begin{center}
		\subfloat[]{\includegraphics[height=0.255\textwidth,width=0.255\textwidth]{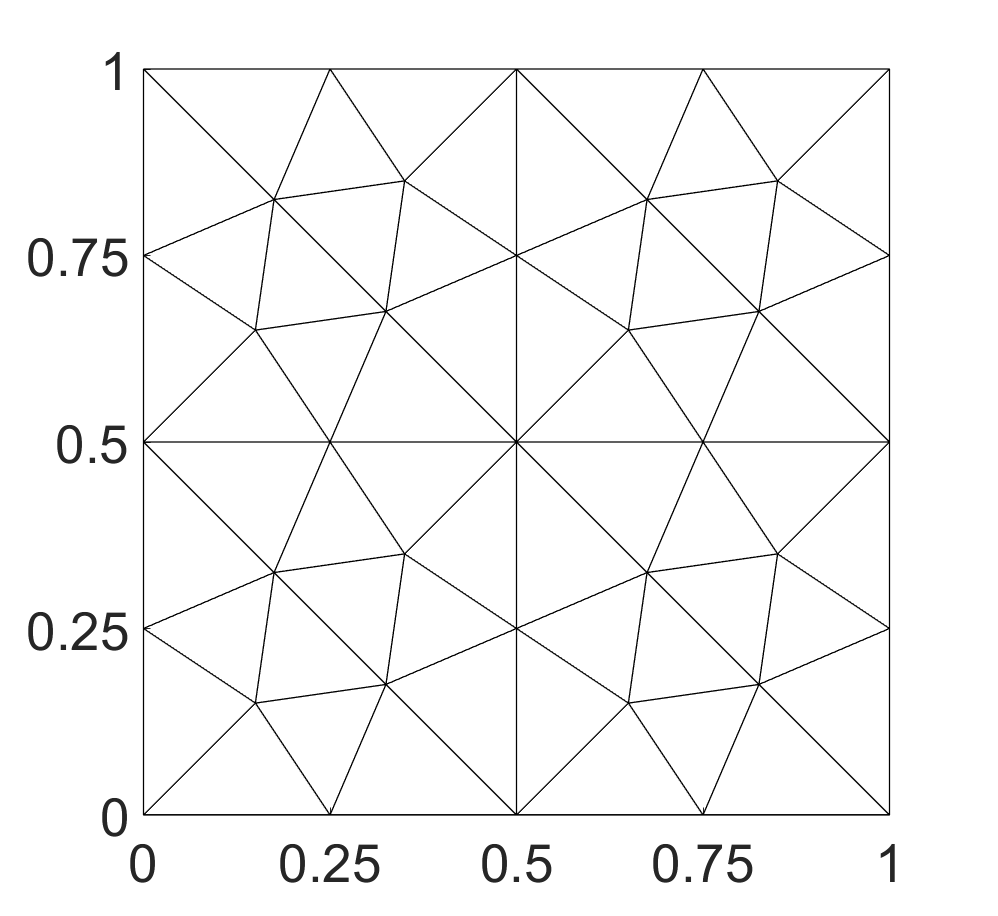}}\hspace{-1em}
		\subfloat[]{\includegraphics[height=0.255\textwidth,width=0.255\textwidth]{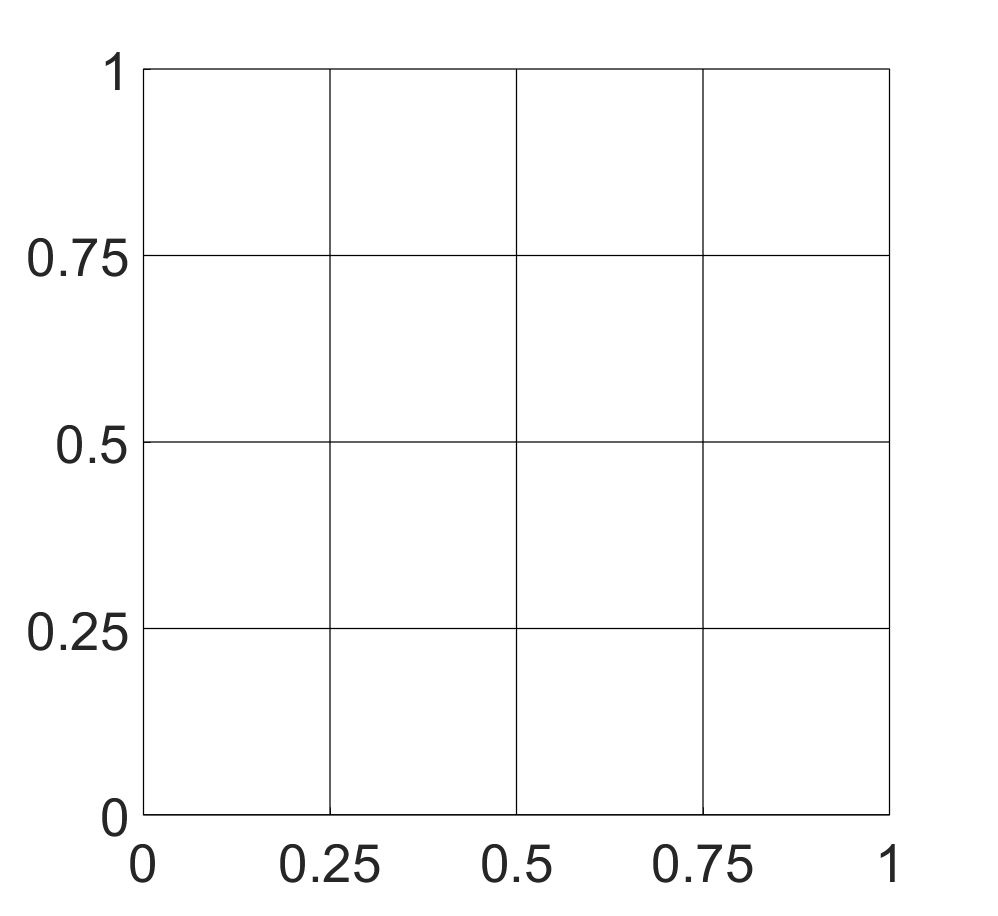}}\hspace{-1em}
		\subfloat[]{\includegraphics[height=0.255\textwidth,width=0.255\textwidth]{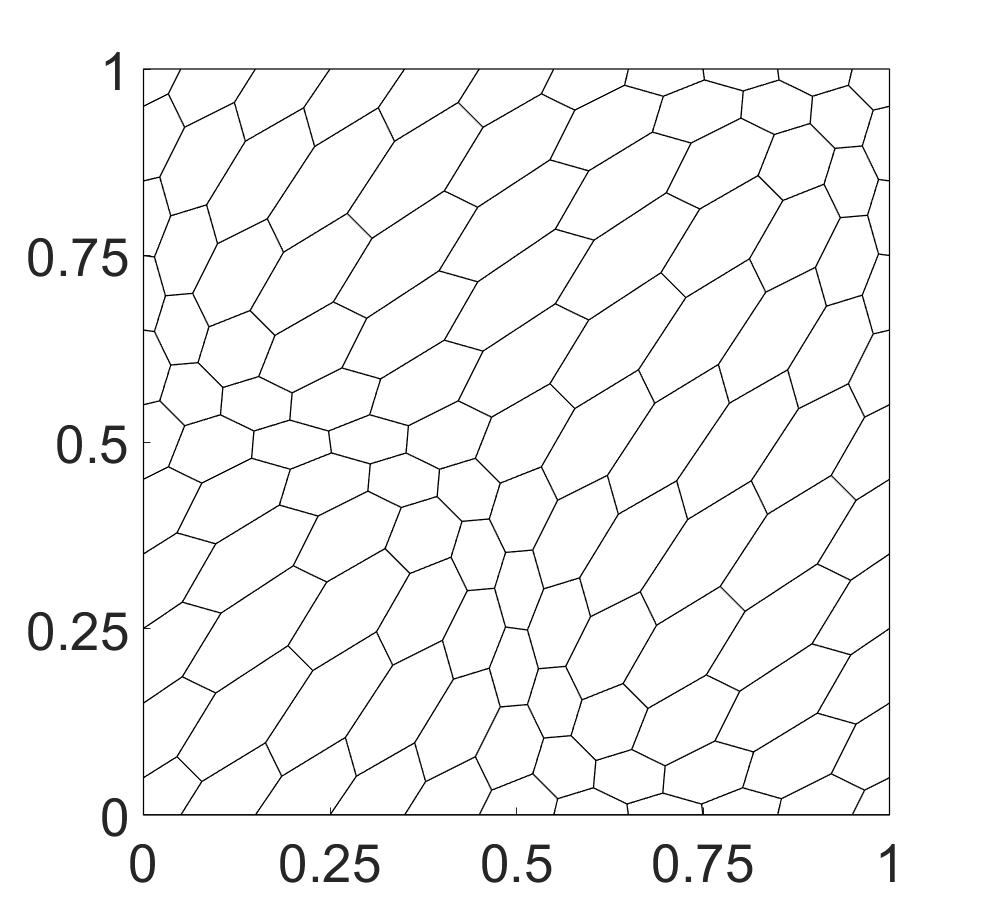}}\hspace{-1em}
		\subfloat[]{\includegraphics[height=0.255\textwidth,width=0.255\textwidth]{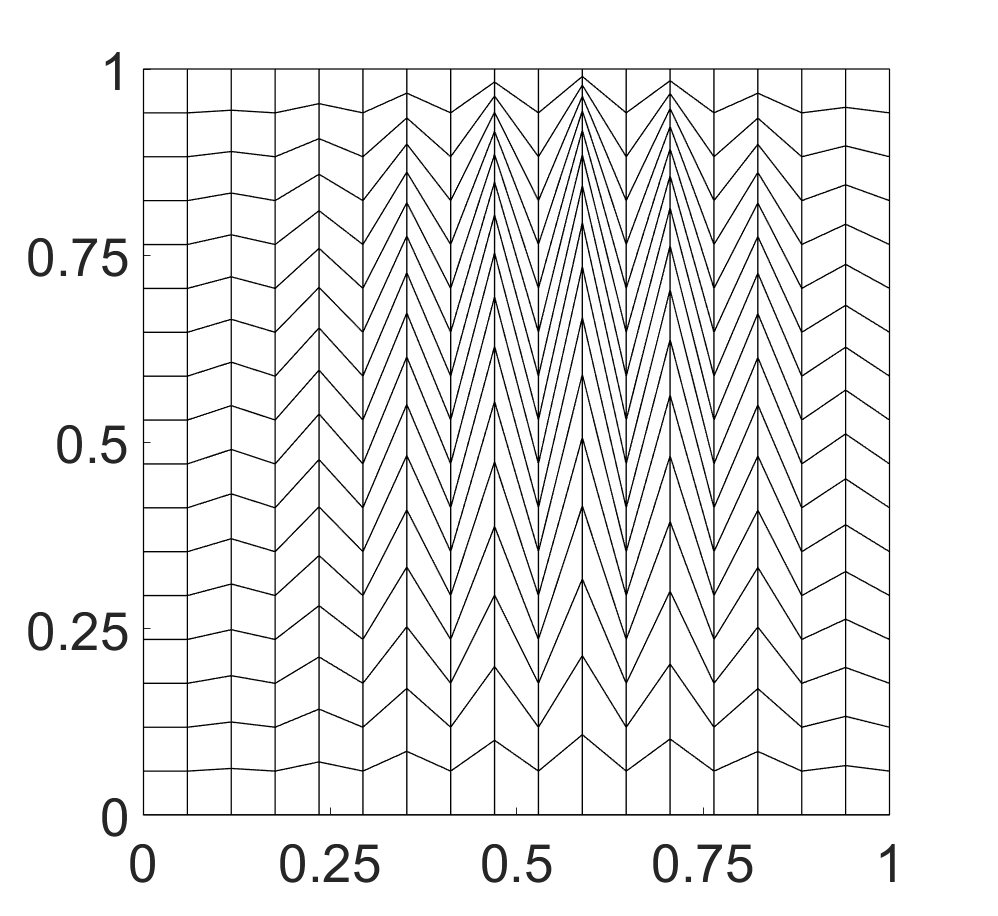}}
		\caption{(a) Triangular, (b) Cartesian, (c) hexagonal and (d) Kershaw initial meshes.}
		\label{fig:Initial_Meshes}
	\end{center}
\end{figure}

\begin{figure}[!ht]
\begin{center}
	\subfloat[]{\includegraphics[height=0.4\textwidth,width=0.5\textwidth]{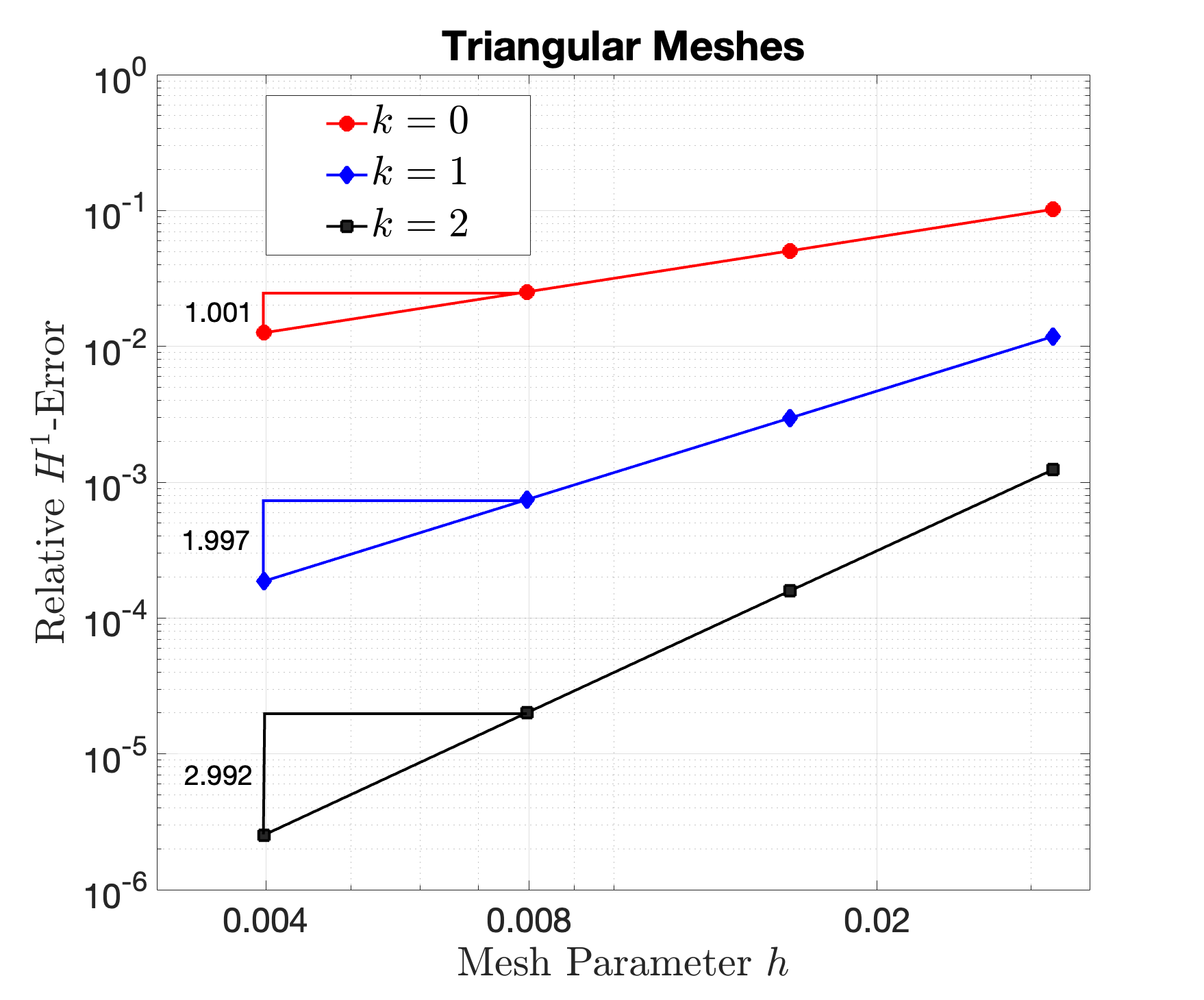}}
	\subfloat[]{\includegraphics[height=0.4\textwidth,width=0.5\textwidth]{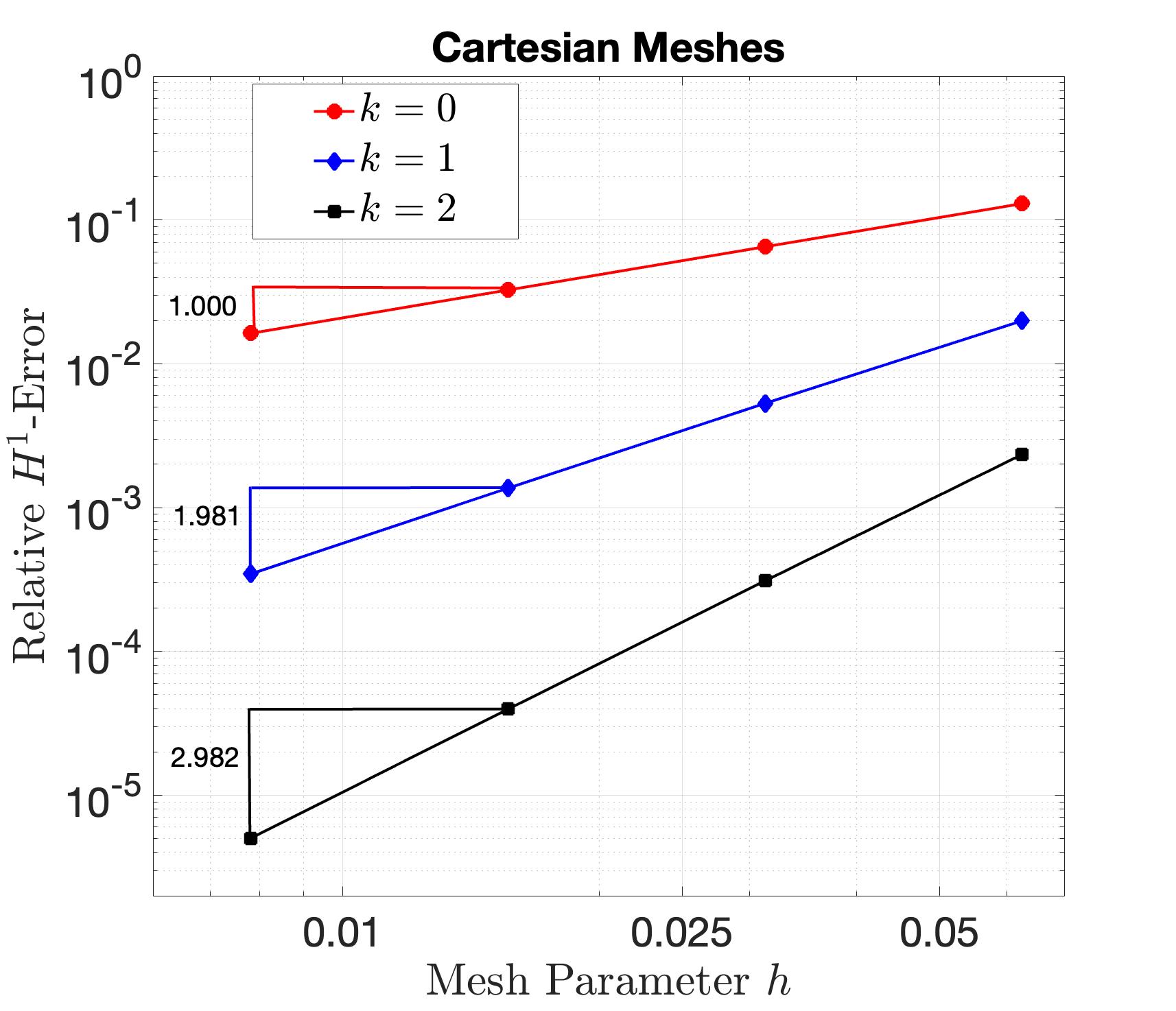}}\\
	\subfloat[]{\includegraphics[height=0.4\textwidth,width=0.5\textwidth]{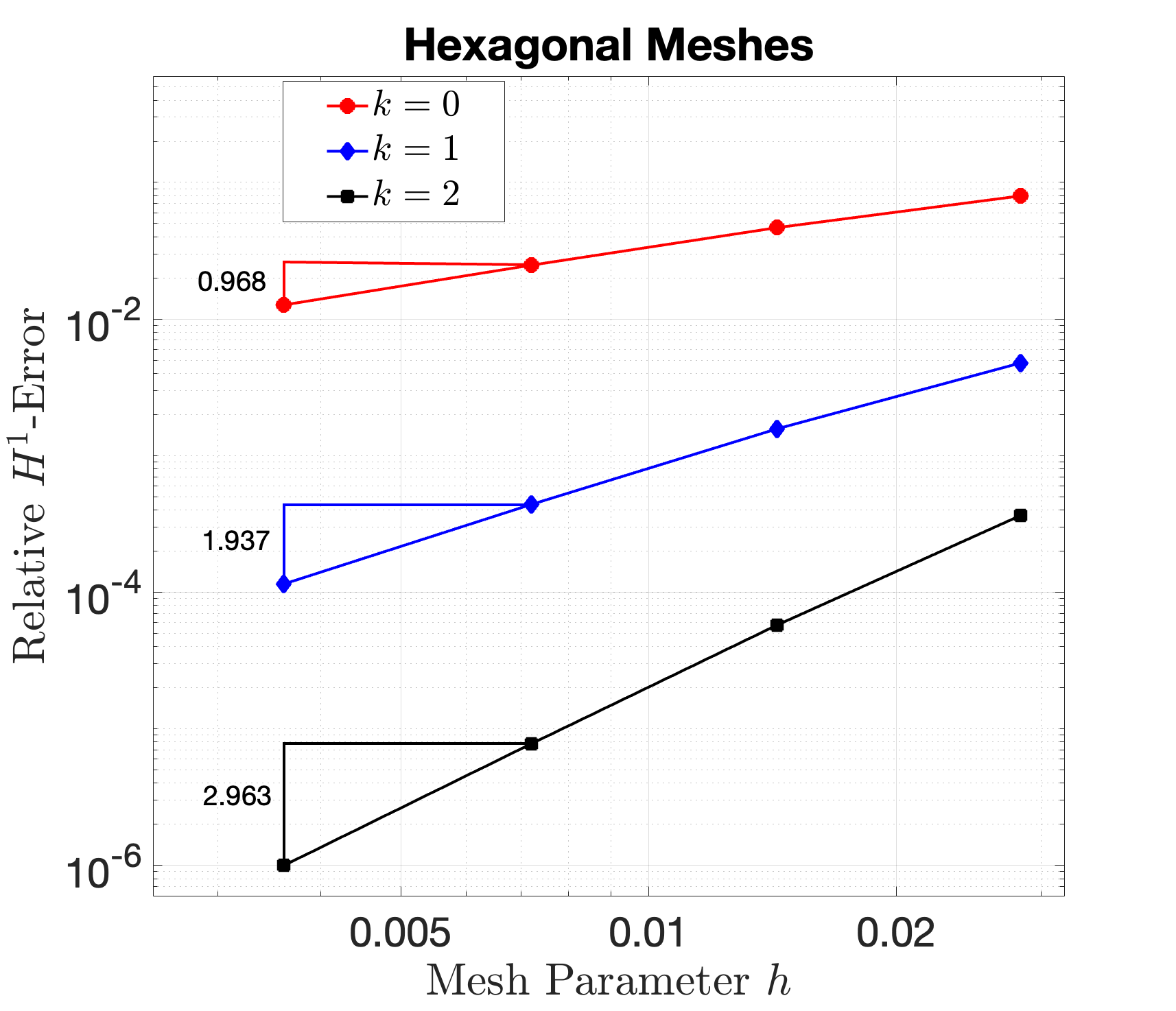}}
	\subfloat[]{\includegraphics[height=0.4\textwidth,width=0.5\textwidth]{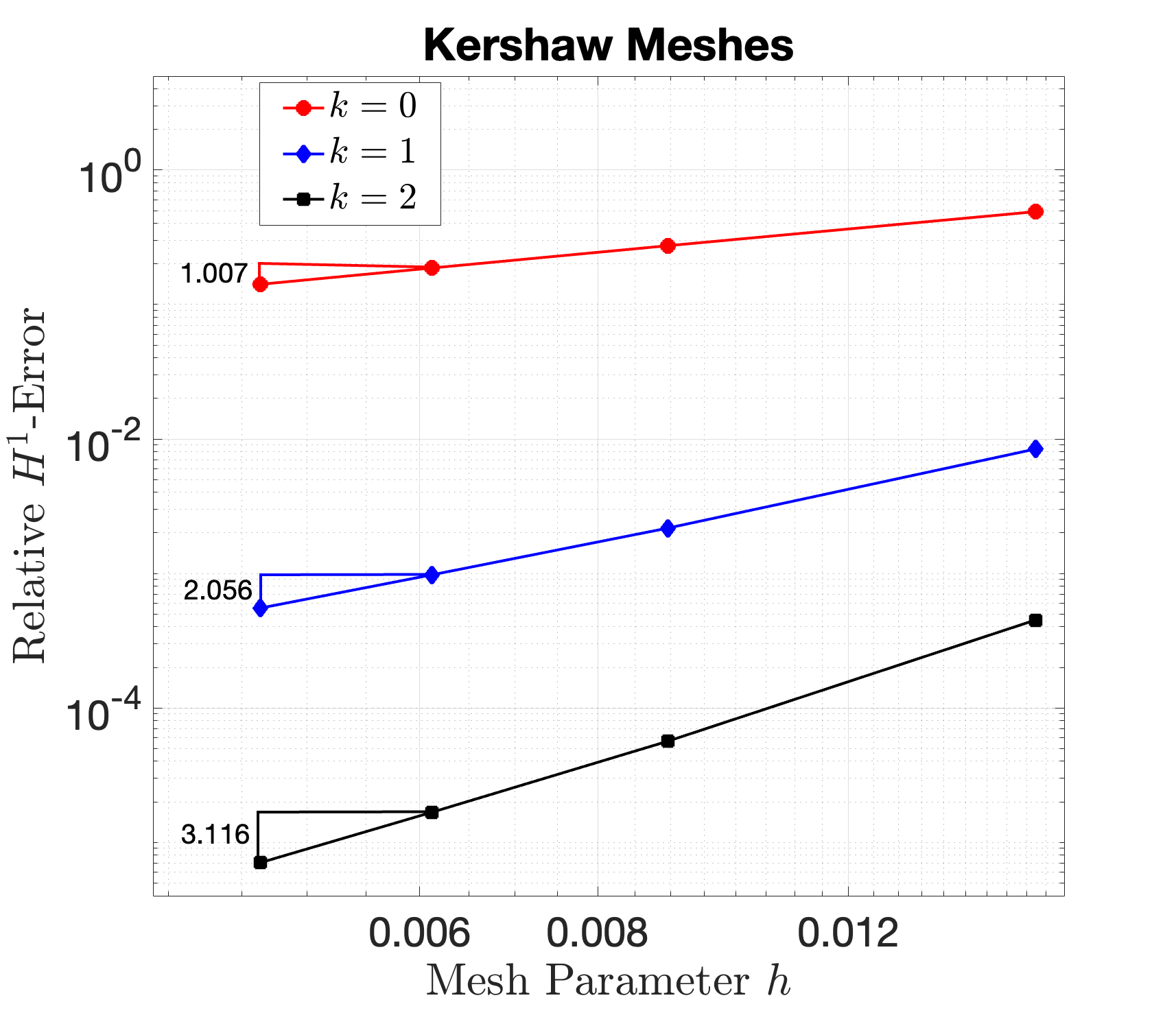}}
	\caption{Convergence histories for the relative discrete $H^1$-error on the (a) triangular, (b)  Cartesian, (c) hexagonal and (d) Kershaw meshes.}
	\label{fig:Kirchhoff_Conv_His}
\end{center}
\end{figure}

\begin{example}
Consider the following nonlocal Kirchhoff type problem: 
 \begin{align}\label{numer_kirchhoff}
    -\left(1+\|\nabla u\|^2\right)\Delta u = f(x)\;\;\text{ in }\; \Omega,\;\; u=0\;\;\text{ in }\;\partial\Omega,
\end{align}
 where $\Omega=(0,1)\times (0,1)\subset \mR^2$ is the unit square domain. Note that with respect to \eqref{cts_kirchhoff}, here $M(d)=1+d$ with $d=\|\nabla u\|^2$ and $L_f=0$. The existence of a unique solution of the continuous problem follows from Remark~\ref{rem_exit_Browder}. The source term $f$ is chosen in such a way that the exact solution reads $u(x,y)=x(1-x)y(1-y)$. We perform numerical tests on sequences of triangular, Cartesian, hexagonal and Kershaw meshes, whose initial meshes are shown in Figure~\ref{fig:Initial_Meshes}. The triangular, Cartesian and Kershaw mesh families correspond to mesh families 1, 2 and 4.1, respectively, of the FVCA5 benchmark \cite{Her_Hub_FVCA_mesh_08}. The (predominantly) hexagonal mesh family originates from the work of \cite{Piet_Lem_HexaMesh_15}. We follow some of the HHO implementation strategies described in \cite{Piet_Jero_HHO_Book_20,TG_GM_TP_22_HHO_Quasi,GM_TG_24_HHO_Strong_Nonlin}. he source code and datasets generated and analyzed for this article are publicly available on the GitHub repository: https://github.com/gourangamallik/HHO-Nonlocal-Kirchhoff. We apply Newton's iteration \eqref{Newton_iter} to compute the discrete solution on a mesh $\cT_h$. Choose the initial guess $\underline{u}_h^{(0)}$ as the HHO solution of the Poisson problem: $-\Delta u=f$ in $\Omega$ and $u=0$ on $\partial\Omega$ with the same load $f$ as in \eqref{numer_kirchhoff}. Set the initial guess for $d^{(0)}=\|\nabla R_h^{k+1}\underline{u}_h^{(0)}\|^2$. The stopping criterion is based on the relative error between two successive iterates, that is,
 \begin{equation*}
 \|\underline{u}_h^{(n)}-\underline{u}_h^{(n-1)}\|_{a,h}/\|\underline{u}_h^{(n)}\|_{a,h}\leq 10^{-10}\quad\text{ for } n=1,2,3,\ldots.
 \end{equation*}
 Under this relative error criterion, the iterative scheme consistently converges in at most five iterations for all four different mesh types.

 In Table~\ref{table:Kirchhoff_Triag}-\ref{table:Kirchhoff_Kers}, we report the relative $H^1$-error $e_h$ and the convergence rate for the triangular, Cartesian, hexagonal and Kershaw mesh families. The corresponding convergence histories of the error $e_h$ with respect to the meshsize $h$ are shown in Figure~\ref{fig:Kirchhoff_Conv_His} for polynomial degrees $k=0,1,2$.
The results show that the empirical convergence rates are approximately $1,2$ and $3$ for $k=0,1$ and $2$, respectively, in all mesh families. The numerical results substantiate the theoretical results of Theorem~\ref{thm:Kirchhoof_err_est}, thus confirming the robustness of the HHO method.
\end{example} 

\begin{table}[!ht]
\caption{Relative $H^1$ errors and convergence rates on the triangular meshes.}\label{table:Kirchhoff_Triag}
\begin{center}
	\begin{tabular}{  c c  c c c c c }     
		\hline
		\multirow{2}{*}{$h$}  & \multicolumn{2}{c}{$k=0$} &\multicolumn{2}{c}{$k=1$} &\multicolumn{2}{c}{$k=2$}\\
		\cline{2-7}
		&$e_h$ & \texttt{rate}  & $e_h$ & \texttt{rate} & $e_h$ & \texttt{rate}\\ 
		\hline %\hline
        %0.0732   & 2.1448e-01   & --     & 1.0481e-01  & --     & 1.9240e-02 & --\\  %\hline 
		0.0318   & 1.0212e--01  & --  & 1.1834e--02 & --  & 1.2383e--03 & -- \\  %\hline 
		0.0159   & 5.0460e--02  & 1.017  & 2.9675e--03 & 1.996  & 1.5857e--04 & 2.965 \\ %\hline
		0.0080   & 2.5166e--02  & 1.004  & 7.4427e--04 & 1.995  & 2.0057e--05 & 2.983 \\ %\hline
		0.0040   & 1.2575e--02  & 1.001  & 1.8643e--04 & 1.997  & 2.5218e--06 & 2.992 \\ \hline
	\end{tabular}
\end{center}
\end{table}

\begin{table}[!ht]
\caption{Relative $H^1$ errors and convergence rates on the Cartesian meshes.}\label{table:Kirchhoff_Cart}
\begin{center}
	\begin{tabular}{  c c  c c c c c }     
		\hline
		\multirow{2}{*}{$h$}  & \multicolumn{2}{c}{$k=0$} &\multicolumn{2}{c}{$k=1$} &\multicolumn{2}{c}{$k=2$}\\
		\cline{2-7}
		&$e_h$ & \texttt{rate}  & $e_h$ & \texttt{rate} & $e_h$ & \texttt{rate}\\ 
		\hline %\hline
        %0.1250   & 2.5582e-01   & --     & 6.7401e-02  & --     & 1.6341e-02 & --\\  %\hline 
		0.0625   & 1.2988e--01  & --  & 1.9894e--02 & --  & 2.3401e--03 & -- \\  %\hline 
		0.0313   & 6.5098e--02  & 0.996  & 5.3052e--03 & 1.907  & 3.0948e--04 & 2.919 \\ %\hline
		0.0156   & 3.2567e--02  & 0.999  & 1.3647e--03 & 1.959  & 3.9701e--05 & 2.963 \\ %\hline
		0.0078   & 1.6285e--02  & 1.000  & 3.4580e--04 & 1.981  & 5.0248e--06 & 2.982 \\ \hline
	\end{tabular}
\end{center}
\end{table}

\begin{table}[!ht]
\caption{Relative $H^1$ errors and convergence rates on the hexagonal meshes.}\label{table:Kirchhoff_Hexa}
\begin{center}
	\begin{tabular}{  c c  c c c c c }     
		\hline
		\multirow{2}{*}{$h$}  & \multicolumn{2}{c}{$k=0$} &\multicolumn{2}{c}{$k=1$} &\multicolumn{2}{c}{$k=2$}\\
		\cline{2-7}
		&$e_h$ & \texttt{rate}  & $e_h$ & \texttt{rate} & $e_h$ & \texttt{rate}\\ 
		\hline %\hline
		0.0283   & 7.9632e--02  & --     & 4.7636e--03 & --     & 3.6400e--04 & --\\  %\hline 
		0.0143   & 4.6656e--02  & 0.785  & 1.5704e--03 & 1.629  & 5.7435e--05 & 2.711 \\ %\hline
		0.0072   & 2.4788e--02  & 0.920  & 4.3927e--04 & 1.854  & 7.7763e--06 & 2.910 \\ %\hline
		0.0036   & 1.2681e--02  & 0.968  & 1.1490e--04 & 1.937  & 9.9923e--07 & 2.963 \\ \hline
	\end{tabular}
\end{center}
\end{table}

\begin{table}[!ht]
\caption{Relative $H^1$ errors and convergence rates on the Kershaw meshes.}\label{table:Kirchhoff_Kers}
\begin{center}
	\begin{tabular}{  c c  c c c c c }     
		\hline
		\multirow{2}{*}{$h$}  & \multicolumn{2}{c}{$k=0$} &\multicolumn{2}{c}{$k=1$} &\multicolumn{2}{c}{$k=2$}\\
		\cline{2-7}
		&$e_h$ & \texttt{rate}  & $e_h$ & \texttt{rate} & $e_h$ & \texttt{rate}\\ 
		\hline %\hline
		0.0162   & 4.9034e--01  & --     & 8.4119e--03 & --     & 4.4978e--04 & --\\  %\hline 
	    0.0090   & 2.7388e--01  & 0.980  & 2.1705e--03 & 2.280  & 5.6402e--05 & 3.494 \\ %\hline
		0.0061   & 1.8697e--01  & 1.002  & 9.7559e--04 & 2.099  & 1.6711e--05 & 3.192 \\ %\hline
		0.0046   & 1.4144e--01  & 1.007  & 5.5181e--04 & 2.056  & 7.0458e--06 & 3.116 \\ \hline
	\end{tabular}
\end{center}
\end{table}

\section{Conclusion}
This article introduces a hybrid high-order finite element approximation for nonlocal nonlinear problems of Kirchhoff type. The existence and uniqueness of the discrete solution are established without imposing additional assumptions on the underlying continuous problem. Moreover, we show that, for a particular type of nonlocal Kirchhoff problem, the existence and uniqueness of continuous and discrete solutions can be established without imposing Assumptions (A1)--(A3). The nonconformity and inconsistency inherent in the HHO method are effectively controlled through a suitable choice of the mesh parameter. Furthermore, we derive error estimates in an $H^1$-like discrete norm for arbitrary-order polynomial approximations. Numerical experiments are performed on four different types of polygonal meshes that confirm the robustness and optimal convergence of the proposed method.

\medskip

\noindent\textbf{Declaration of competing interest}\\
The author declares that he has no conflict of interest.

\medskip

\noindent\textbf{Acknowledgements}\\
The author acknowledges the DST-SERB-MATRICS grant MTR/2023/000681 for financial support and the DST FIST grant SR/FST/MS II/2023/139 for partial financial support.

\bigskip

\bibliographystyle{siam}
\bibliography{MyBiblio,HHO_Bib}
%\bibliography{Kirchhoff}

\appendix
\section{Continuity result of $\mathcal{S}$:}
In this appendix, we provide the proof of the continuity of the nonlinear map $\mathcal{S}:\underline{U}_{h,0}^k\to \underline{U}_{h,0}^k$ as defined in \eqref{defn_S}. Recall the definition of $\mathcal{S}(\underline{w}_h)\in \underline{U}_{h,0}^k$ from \eqref{defn_S}
\begin{equation*}
    a_h(\mathcal{S}(\underline{w}_h),\underline{v}_h)=Q_{w}(\underline{v}_h)\quad\forall\underline{v}_h\in \underline{U}_{h,0}^k,
\end{equation*} 
% From the boundedness of $Q_w$ in \eqref{Q_bdd}, we have 
% \begin{align}
%     \|\mathcal{S}(\underline{w}_h)\|_{a,h}&=\sup_{\|\underline{v}_h\|_{a,h}\leq 1}a_h(\mathcal{S}(\underline{w}_h),\underline{v}_h)=\sup_{\|\underline{v}_h\|_{a,h}\leq 1}Q_{w}(\underline{v}_h)\notag\\
%     &\leq  M\left(\|\nabla R_h^{k+1}\underline{w}_h\|^2\right)\|\underline{w}_h\|_{a,h}+\lambda_{1,h}^{-1/2}a\|w_h\|+\|f(x,0)\|.
% \end{align}
where $Q_{w}(\underline{v}_h)=M\left(\|\nabla R_h^{k+1}\underline{w}_h\|^2\right)a_h(\underline{w}_h,\underline{v}_h)-(f(w_h),v_h)$. Let $\underline{w}_h\in \underline{U}_{h,0}$ be any discrete function. Let $\epsilon>0$ be given. 
%From the definition of $Q_w$ in \eqref{defn_Qw} and $\mathcal{S}(\underline{w}_h)$ in \eqref{defn_Rw}, we have
Now
\begin{align}\label{diff_R}
    a_h(\mathcal{S}(\underline{w}_h)-\mathcal{S}(\underline{\psi}_h),\underline{v}_h)&=M\left(\|\nabla R_h^{k+1}\underline{w}_h\|^2\right)a_h(\underline{w}_h,\underline{v}_h)-M\left(\|\nabla R_h^{k+1}\underline{\psi}_h\|^2\right)a_h(\underline{\psi}_h,\underline{v}_h)\notag\\
    &\qquad+(f(w_h)-f(\psi_h),v_h).
\end{align}
If $\underline{w}_h=\underline{0}$, then
\begin{align*}
    |a_h(\mathcal{S}(\underline{0}_h)-\mathcal{S}(\underline{\psi}_h),\underline{v}_h)|&=|M\left(\|\nabla R_h^{k+1}\underline{\psi}_h\|^2\right)a_h(\underline{\psi}_h,\underline{v}_h)|+|(f(0)-f(\psi_h),v_h)|\notag\\
    &\leq M\left(\|\nabla R_h^{k+1}\underline{\psi}_h\|^2\right)\|\underline{\psi}_h\|_{a,h}\|\underline{v}_h\|_{a,h}+L_f\|0-\psi_h\|\|v_h\|\notag\\
    &\leq \left(M\left(\|\nabla R_h^{k+1}\underline{\psi}_h\|^2\right)+\frac{L_f}{\lambda_{1,h}}\right)\|\underline{\psi}_h\|_{a,h}\|\underline{v}_h\|_{a,h}.
\end{align*}
Since $M$ is continuous, there exists $\delta>0$ such that $\left(M\left(\|\nabla R_h^{k+1}\underline{\psi}_h\|^2\right)+\frac{L_f}{\lambda_{1,h}}\right)\|\underline{\psi}_h\|_{a,h}<\epsilon$ whenever $\|\underline{0}-\underline{\psi}_h\|_{a,h}<\delta$. Consequently, $\|\mathcal{S}(\underline{0}_h)-\mathcal{S}(\underline{\psi}_h)\|_{a,h}<\epsilon$ whenever $\|\underline{0}-\underline{\psi}_h\|_{a,h}<\delta$. Therefore, $\mathcal{S}$ is continuous at $\underline{0}\in \underline{U}_{h,0}^k$. Now, for $\underline{w}_h\neq\underline{0}$, we have
\begin{align}\label{M_diff}
 &M\left(\|\nabla R_h^{k+1}\underline{w}_h\|^2\right)a_h(\underline{w}_h,\underline{v}_h)-M\left(\|\nabla R_h^{k+1}\underline{\psi}_h\|^2\right)a_h(\underline{\psi}_h,\underline{v}_h)\notag\\
 &= \left(M\left(\|\nabla R_h^{k+1}\underline{w}_h\|^2\right)-M\left(\|\nabla R_h^{k+1}\underline{\psi}_h\|^2\right)\right)a_h(\underline{w}_h,\underline{v}_h)+M\left(\|\nabla R_h^{k+1}\underline{\psi}_h\|^2\right)a_h(\underline{w}_h-\underline{\psi}_h,\underline{v}_h).
\end{align}
Since $M$ is continuous, there exists $\delta_1>0$ such that
\begin{align*}
    \left|M\left(\|\nabla R_h^{k+1}\underline{w}_h\|^2\right)-M\left(\|\nabla R_h^{k+1}\underline{\psi}_h\|^2\right)\right|<\frac{\epsilon}{2\|\underline{w}\|_{a,h}},\; \text{ whenever } \|\underline{w}_h-\underline{\psi}_h\|_{a,h}<\delta_1.
\end{align*}
Using the boundedness of $a(\cdot,\cdot)$, we obtain the following from \eqref{M_diff} as 
\begin{align*}
&\left|M\left(\|\nabla R_h^{k+1}\underline{w}_h\|^2\right)a_h(\underline{w}_h,\underline{v}_h)-M\left(\|\nabla R_h^{k+1}\underline{\psi}_h\|^2\right)a_h(\underline{\psi}_h,\underline{v}_h)\right|\\
 &\leq\left(\frac{\epsilon}{2}+M(\|\nabla R_h^{k+1}\underline{w}_h\|^2)\|\underline{w}_h-\underline{\psi}_h\|_{a,h}\right)\|\underline{v}_h\|_{a,h}.
\end{align*}
Moreover, using the Lipschitz continuity of $f$, we have
\begin{equation*}
    |(f(w_h)-f(\psi_h),v_h)|\leq L_f \|w_h-\psi_h\|\|v_h\|\leq \frac{L_f}{\lambda_{1,h}} \|\underline{w}_h-\underline{\psi}_h\|_{a,h}\|\underline{v}_h\|_{a,h}.
\end{equation*}
Finally, combining the above estimates in \eqref{diff_R} and using the boundedness of $a_h(\cdot,\cdot)$, we have 
\begin{equation*}
    \|\mathcal{S}(\underline{w}_h)-\mathcal{S}(\underline{\psi}_h)\|_{a,h}\leq \frac{\epsilon}{2}+M(\|\nabla R_h^{k+1}\underline{w}_h\|^2)\|\underline{w}_h-\underline{\psi}_h\|_{a,h}+\frac{L_f}{\lambda_{1,h}}\|\underline{w}_h-\underline{\psi}_h\|_{a,h}.
\end{equation*}
Therefore, there exists $\delta_2>0$ such that $M(\|\nabla R_h^{k+1}\underline{w}_h\|^2)\|\underline{w}_h-\underline{\psi}_h\|_{a,h}+\frac{a}{\lambda_{1,h}}\|\underline{w}_h-\underline{\psi}_h\|_{a,h}<\frac{\epsilon}{2}$ whenever $\|\underline{w}_h-\underline{\psi}_h\|_{a,h}<\delta_2$. Choose $\delta=\min \{\delta_1,\delta_2\}$. Then
\begin{equation*}
    \|\mathcal{S}(\underline{w}_h)-\mathcal{S}(\underline{\psi}_h)\|_{a,h}<\epsilon\;\text{ whenever } \|\underline{w}_h-\underline{\psi}_h\|_{a,h}< \delta.
\end{equation*}
This implies that $\mathcal{S}$ is continuous on $\underline{U}_{h,0}^k$.
\end{document}